\newif\ifpreprint
\preprinttrue

\ifpreprint
\documentclass[12pt]{article}
\usepackage{fullpage}
\else
\RequirePackage{fix-cm}
\documentclass[smallextended]{svjour3}       \smartqed                                                      \journalname{Mathematical Programming Computation}
\fi

\usepackage{graphicx}
\usepackage{pgfplots}
\pgfplotsset{compat=1.17}
\usepgfplotslibrary{groupplots}
\usepackage{tikz}
\usetikzlibrary{arrows.meta, positioning, fit, backgrounds, calc}
\definecolor{fomtoolsblue}{RGB}{26,95,150}
\definecolor{cbBlue}{RGB}{0,114,178}     \definecolor{cbOrange}{RGB}{230,159,0}   \definecolor{cbPurple}{RGB}{204,121,167} 

\ifpreprint \usepackage{amsthm}

\fi

\usepackage{amsmath,amssymb,enumitem,mathtools}
\renewcommand\arraystretch{1.2}
\usepackage{multirow}
\usepackage{algorithm}\usepackage[noend]{algpseudocode} 

\usepackage{listings}
\definecolor{codekeyword}{HTML}{008000}
\definecolor{codecomment}{HTML}{3D7B7B}
\definecolor{codestring}{HTML}{BA2121}
\definecolor{codenumber}{HTML}{666666}
\lstdefinestyle{fomcode}{
	basicstyle=\small\ttfamily,
	columns=fullflexible,
	keepspaces=true,
	upquote=true,
	showstringspaces=false,
	keywordstyle=\bfseries\color{codekeyword},
	keywordstyle=[2]\color{codekeyword},
	commentstyle=\itshape\color{codecomment},
	stringstyle=\color{codestring},
	numberstyle=\tiny\color{codenumber},
	numbersep=8pt,
	frame=leftline,
	framesep=2mm,
	xleftmargin=17pt,
}
\lstdefinelanguage{fomPython}[]{Python}{
	morekeywords=[2]{range, None, True, False},
}
\makeatletter
\newcommand{\codelisting}[2][]{\par\addvspace{\medskipamount}\noindent\begin{minipage}{\linewidth}\lstinputlisting[style=fomcode, aboveskip=0pt, belowskip=0pt, #1]{#2}\end{minipage}\par\addvspace{\medskipamount}\@doendpe}
\makeatother

\usepackage[legacy]{csvsimple}	\usepackage{booktabs}   \usepackage{array}      \usepackage{adjustbox}  

\usepackage[colorlinks=true, citecolor=cyan, linkcolor=cyan, urlcolor=cyan]{hyperref}

\usepackage[xindy, acronyms]{glossaries}   

\newif\ifshownotes
\shownotesfalse \ifshownotes
	\newcommand{\vnote}[1]{\textcolor{magenta}{\textbf{[Vinit: #1]}}}
	\newcommand{\bnote}[1]{\textcolor{cyan}{\textbf{[Bartolomeo: #1]}}}
	\newcommand{\dnote}[1]{\textcolor{purple}{\textbf{[Daniel: #1]}}}
	\newcommand{\jnote}[1]{\textcolor{orange}{\textbf{[Jisun: #1]}}}
\else
	\newcommand{\vnote}[1]{\unskip}
	\newcommand{\bnote}[1]{\unskip}
	\newcommand{\dnote}[1]{\unskip}
	\newcommand{\jnote}[1]{\unskip}
\fi

\newcommand{\eg}{{\it e.g.}}
\newcommand{\ie}{{\it i.e.}}

\newcommand{\ones}{\mathbf 1}
\newcommand{\reals}{{\mathbf R}}

\newcommand{\symm}{{\mathbf S}}  \newcommand{\NPhard}{\mbox{$\mathcal{NP}$-hard}}  \newcommand{\prob}{\mathop{\mathbf{prob}}}

\newcommand{\tpose}{T}

\newcommand{\diag}{\mathop{\bf diag}}
\newcommand{\lambdamax}{{\lambda_{\rm max}}}

\newcommand{\argmin}{\mathop{\rm argmin}}

\newcommand{\prox}{{\bf prox}}

\newcommand{\sign}{\mathop{\bf sign}}

\newcommand{\softthresh}{\mathcal{T}}
\newcommand{\relu}{\mathrm{relu}}
\newcommand{\sat}{\mathrm{sat}}

\newcommand{\Sec}{Section\;}

\ifpreprint

\theoremstyle{definition}

\theoremstyle{remark}

\fi

\newacronym{LP}{LP}{linear program}
\newacronym{QP}{QO}{quadratic optimization}
\newacronym{MIQP}{MIQP}{mixed-integer quadratic program}
\newacronym{MIP}{MIP}{mixed-integer program}
\newacronym{SDP}{SDP}{semidefinite program}
\newacronym{MILP}{MILP}{mixed-integer linear program}
\newacronym{MINLP}{MINLP}{mixed-integer nonlinear program}
\newacronym{sBB}{sBB}{spacial branch and bound}
\newacronym{PWA}{PWA}{piecewise affine}
\newacronym{PDHG}{PDHG}{primal-dual hybrid gradient}
\newacronym{SVM}{SVM}{support vector machines}
\newacronym{ReLU}{ReLU}{rectified linear unit}
\newacronym{CPU}{CPU}{central processing unit}
\newacronym{GPU}{GPU}{graphics processing unit}
\newacronym{MPC}{MPC}{model predictive control}
\newacronym{ADMM}{ADMM}{alternating direction method of multipliers}
\newacronym{ADP}{ADP}{approximate dynamic programming}
\newacronym{FPGA}{FPGA}{field-programmable gate array}
\newacronym{ABI}{ABI}{application binary interface}
\newacronym{IR}{IR}{intermediate representation}
\newacronym{LLM}{LLM}{large language model}
\newacronym{CG}{CG}{conjugate gradient}
\newacronym{CSR}{CSR}{compressed sparse row}
\newacronym{CSC}{CSC}{compressed sparse column}
\newacronym{KKT}{KKT}{Karush--Kuhn--Tucker}
\newacronym{RCM}{RCM}{reverse Cuthill--McKee}

\newcommand{\ldlt}{\mathrm{LDL}^\tpose}

\newcommand{\fomtools}{\texttt{FOMTools}}

\newcommand{\paperauthors}{Vinit~Ranjan\textsuperscript{*},
    Daniel~Deza\textsuperscript{\dag},
    Anushka~Acharya\textsuperscript{\ddag},
    Jisun~Park\textsuperscript{\dag},
    Bartolomeo~Stellato\textsuperscript{\dag}}
\newcommand{\paperaffiliations}{\textsuperscript{*}Sloan School of Management, Massachusetts Institute of Technology\\
    \textsuperscript{\dag}Department of Operations Research and Financial Engineering, Princeton University\\
    \textsuperscript{\ddag}Data and Intelligent Systems, Princeton University\\[0.5ex]
    {\fontencoding{T1}\ttfamily
    vranjan1@mit.edu\\
    \{dd7022, aa9078, jisunpark, bstellato\}@princeton.edu}}

\title{ FOMTools: Development and AI-Assisted Code Generation for First-Order Methods}
\ifpreprint
\author{\begin{minipage}{\dimexpr\textwidth-2\tabcolsep\relax}
        \centering\normalsize
        \paperauthors\\[0.75em]
        \small\paperaffiliations
    \end{minipage}}
\else
\author{\paperauthors}
\titlerunning{FOMTools: Development and Code Generation}
\authorrunning{Vinit Ranjan et al.}
\fi

\begin{document}

\maketitle

\ifpreprint
\let\oldtabular\tabular
\let\endoldtabular\endtabular
\renewenvironment{tabular}{\renewcommand\arraystretch{1.0}\oldtabular}{\endoldtabular}
\fi

\begin{abstract}
We present FOMTools, an open-source Python package for developing, analyzing, and deploying first-order algorithms for parametric linear and convex quadratic optimization.
The user specifies an iteration in Python, with syntax close to its mathematical description.
The package represents the iteration as a computation graph of affine and piecewise affine operations.
To analyze worst-case performance, FOMTools formulates a mixed-integer verification problem that maximizes a convergence residual over a specified set of problem instances and initial iterates.
For probabilistic performance analysis, a batched GPU simulator collects algorithm trajectories and computes statistical bounds on the probability of exceeding a specified convergence tolerance.
From the same algorithm description, FOMTools generates specialized solver code for embedded hardware, multicore CPUs, and GPUs.
The package first generates the numerical data, factorizations, and matrix-vector routines without any agent.
A coding agent then translates the iteration into C or CUDA code and iteratively optimizes the implementation.
After each phase, the package runs the generated code on sampled inputs and compares it with a reference execution.
FOMTools also serves as a CVXPYgen backend, generating and optimizing solvers for parametric CVXPY problems.
We demonstrate FOMTools on sparse coding, network optimization, optimal control, and portfolio rebalancing.
The verification bounds reveal worst-case behavior that sampling misses.
The generated solvers run up to two orders of magnitude faster than NumPy implementations, and they outperform hand-written GPU code and the OSQP solvers generated by CVXPYgen.
\end{abstract}

\section{Introduction}\label{sec:introduction}

First-order methods are widely used for large-scale convex optimization because their iterations often require only matrix-vector products and simple projections or proximal operations.
They are also useful when an application repeatedly solves problems with the same structure, as in model predictive control~\cite{borrelliPredictiveControlLinear2017}, sparse signal reconstruction~\cite{tibshiraniRegressionShrinkageSelection1996a,xuDOAEstimationBased2012}, network optimization~\cite{bertsekas1991linear}, and portfolio rebalancing~\cite{boydMultiPeriodTradingConvex2017}.
We call such a collection a \emph{parametric family}, with the varying data restricted to a specified parameter set.
We focus on linear and convex quadratic programs with fixed matrices and varying objective and constraint vectors.
A solution from one instance can initialize the next solve, often reducing the number of iterations through a warm start.

An application may limit the solver to a fixed iteration budget to meet a runtime or throughput requirement.
For example, real-time electricity markets require repeated optimal power flow solves at five-minute intervals~\cite{vanHentenryckMachineLearningOptimal2021}.
The MIT Cheetah~3 robot controller solves its predictive-control quadratic programs in under $1$ millisecond within a $20$--$30$ Hz planning loop~\cite{diCarloDynamicLocomotion2018}.
The application's tolerance for failure determines the required accuracy guarantee within this budget.
We seek either a bound that holds uniformly over the family or a bound on the probability of exceeding an accuracy tolerance under a specified distribution of instances and initial iterates.

Classical convergence analysis bounds performance after a finite number of iterations over function classes, such as smooth convex functions~\cite{beckFastIterativeShrinkageThresholding2009a}.
The performance estimation problem (PEP) framework~\cite{droriPerformanceFirstorderMethods2014,taylorExactWorstCasePerformance2017,taylorSmoothStronglyConvex2017} computes worst-case bounds over such classes through semidefinite optimization.
These bounds need not capture a particular family's fixed matrices, parameter restrictions, or warm starts, and can therefore be pessimistic for that family~\cite{ranjanVerificationFirstorderMethods2025}.
In contrast, the verification problem~\cite{ranjanVerificationFirstorderMethods2025,ranjanExactVerificationFirstOrder2025} maximizes a performance metric of a specified method over a given parameter set and initial iterate set, after a fixed number of iterations.
Mixed-integer formulations and semidefinite relaxations provide bounds on this worst-case performance.

Constructing such a verification problem requires substantial manual work.
The iterations must be encoded as optimization constraints, together with the admissible parameters and initial iterates.
Obtaining a compact formulation also requires eliminating auxiliary variables and tightening bounds on intermediate iterates~\cite{ranjanExactVerificationFirstOrder2025}.
Adding momentum or changing a regression penalty from squared $\ell_2$ regularization to $\ell_1$ regularization requires revising this formulation.

Deployment requires a second translation of the algorithm.
A Python prototype is convenient for analysis, but embedded hardware targets often need static memory allocation, minimal external library dependencies, and predictable runtime.
A specialized implementation can exploit fixed problem data, sparsity, and reusable factorizations.
Existing code generators automate deployment for many problem families~\cite{banjacEmbeddedCodeGeneration2017,mattingleyCVXGENCodeGenerator2012,schallerEmbeddedCodeGeneration2022}.
Supporting a new iteration or hardware target can still require writing a specialized version of the same algorithm.
Each change must then be reflected in both the verification formulation and the deployed code.
Coding agents based on large language models can now write and tune numerical code from execution and timing feedback~\cite{pressAlgoTuneCanLanguage2025,zimmerAgenticResearcherPractical}.
Their output, however, must be checked against the algorithm it is meant to implement.

In this paper, we present \fomtools{}, an open-source Python package that automates analysis and deployment from one description of a first-order method.
\begin{itemize}
    \item \emph{Method specification.}
    Users specify the algorithm iterations, parameter set, and initial iterate set in Python, with syntax close to the mathematical notation.
    A computation graph records affine and piecewise affine operations, including gradient, momentum, projection, and proximal updates.
    The package substitutes affine expressions into subsequent operations to construct compact verification models and extract the calculations needed by the generated solver.
    \item \emph{Performance analysis.}
    \fomtools{} provides worst-case and probabilistic analysis for a chosen iteration budget and performance metric.
    For worst-case analysis, it constructs a mixed-integer programming (MIP) formulation of the verification problem.
    Bounds on intermediate iterates and convex relaxations of the steps strengthen the formulation to help the solver exclude impossible trajectories.
    For probabilistic analysis, a batched GPU simulator runs the algorithm on inputs sampled from a specified distribution of parameters and initial iterates.
    It computes statistical bounds on the probability of exceeding a prescribed performance tolerance.
    \item \emph{AI-assisted code generation.}
    From the specified algorithm iterations, \fomtools{} generates C code for embedded hardware, OpenMP code for multicore CPUs, and CUDA code for GPUs.
    The package computes factorizations and writes the data and matrix-vector routines needed by the solver.
    An LLM coding agent first translates the iteration into solver code, then profiles and optimizes that code for the target hardware.
    After each phase, the package checks numerical agreement with a reference execution on sampled inputs.
    The package drives the agent through a common interface, so any coding agent can be used.
    An integration with CVXPYgen also generates solvers from problem formulations written in CVXPY.
    \item \emph{Numerical experiments.}
    Experiments on Lasso, network flow, and optimal control compare worst-case and probabilistic bounds, time the generated solvers, and evaluate closed-loop control under a fixed computation budget.
    The generated C solvers are up to $162\times$ faster than NumPy implementations, and code optimization gives up to an additional $12\times$ over the translated solver.
    Generating each of these solvers takes minutes of agent time.
    On a Lasso instance with $40{,}000$ variables, the generated CUDA solvers are $3.3$ to $4.6$ times as fast as hand-written CuPy code.
    On portfolio rebalancing with up to $3000$ assets, solvers generated through \fomtools~are up to $6.2\times$ faster than the OSQP solvers that CVXPYgen generates.
\end{itemize}
The package, together with its documentation and examples, is available at \url{https://github.com/stellatogrp/fomtools}.
The code to reproduce the experiments is available at \url{https://github.com/stellatogrp/fomtools_experiments}.

\subsection{Related work}\label{subsec:related_work}

\paragraph{Performance estimation of first-order methods.}
The PEP framework formulates worst-case performance bounds over a function class as semidefinite programs (SDPs)~\cite{droriPerformanceFirstorderMethods2014,taylorExactWorstCasePerformance2017,taylorSmoothStronglyConvex2017}.
This approach has been used to derive convergence rates and choose coefficients for proximal and operator-splitting methods~\cite{collaAutomaticPerformanceEstimation2023,dasguptaBranchandboundPerformanceEstimation2024,kimDomainSpecificHarnessEndtoEnd2026,parkExactOptimalAccelerated2022,ryuOperatorSplittingPerformance2020,zamaniExactWorstcaseConvergence2024}.
The PESTO~\cite{taylorPerformanceEstimationToolbox2017} and PEPit~\cite{goujaudPEPitComputerassistedWorstcase2024} toolboxes automate these formulations from a user description of the algorithm.
A related SDP-based approach uses integral quadratic constraints (IQCs) to model optimization algorithms as dynamical systems and analyze them using control theory~\cite{lessardAnalysisDesignOptimization2016,taylorLyapunovFunctionsFirstOrder2018}.

\paragraph{Performance verification.}
Verification problems specialize the analysis to a given parametric family and set of initial iterates~\cite{ranjanVerificationFirstorderMethods2025,ranjanExactVerificationFirstOrder2025}.
They retain fixed problem data and use norm-based performance metrics, with warm starts represented by constraints on the initial iterate and its relation to the parameters.
The worst-case analysis in \fomtools{} uses the exact mixed-integer linear programming (MILP) formulations of~\cite{ranjanExactVerificationFirstOrder2025}.
These formulations encode gradient, projection, and proximal iterations through affine and piecewise affine constraints, such as those representing soft-thresholding.

These formulations are closely related to neural network verification, which checks output properties over a specified input set.
Neural network verification combines mixed-integer encodings of rectified linear units~\cite{JMLR:v23:22-0277,fischettiDeepNeuralNetworks2018,tjeng2018evaluating} and convex relaxations~\cite{brownUnifiedViewSDPbased2022} with branch-and-bound algorithms~\cite{bunelUnifiedViewPiecewise2018,shiNeuralNetworkVerification2025}.
GPU-based bound propagation~\cite{wangBetaCROWNEfficientBound2021} and cutting planes~\cite{zhangGeneralCuttingPlanes2022} strengthen these relaxations.
For first-order methods, the computation graph consists of algorithm iterations in place of neural network layers.
When the iteration is averaged or contractive, bounds from operator theory restrict the residuals and distances between iterates.
Combining these bounds with mixed-integer formulations enables verification over tens of iterations in the experiments of~\cite{ranjanExactVerificationFirstOrder2025}.

\paragraph{Data-driven and probabilistic approaches.}
A worst-case bound can be pessimistic for an application if the hardest instances occur rarely.
Sample-based probabilistic analysis uses observed performance to bound the probability of exceeding a tolerance under the sampling distribution.
Combining PEP with distributionally robust optimization extends this analysis to unknown distributions of problem instances~\cite{parkDatadrivenAnalysisFirstOrder2025,ranjanDistributionallyRobustLearningOptimize2026}.
PAC-Bayes generalization bounds provide another way to control the fixed-point residual after a specified number of iterations~\cite{sambharyaLearningWarmStartFixedPoint2024}.
In \fomtools{}, batched GPU execution computes many sample trajectories at once.
The package uses their residuals to compute statistical bounds on the probability of exceeding a prescribed tolerance~\cite{langfordNotBoundingTrue2001,sambharyaDataDrivenPerformanceGuarantees2025}.

\paragraph{Domain-specific languages and software for optimization.}
Modeling languages let users specify optimization problems in a form close to their mathematical description, then convert them to the standard form expected by a solver.
Examples include CVX~\cite{grantCVXMatlabSoftware2014}, CVXPY~\cite{agrawalRewritingSystemConvex2018a,diamondCVXPYPythonEmbeddedModeling2016a}, JuMP~\cite{dunningJuMPModelingLanguage2017}, Convex.jl~\cite{udellConvexOptimizationJulia}, and YALMIP~\cite{lofbergYALMIPToolboxModeling2004}.
\fomtools{} applies this idea to algorithm iterations, which it records as a computation graph.
Its affine expressions build on CVXPY's representation and PEPit's variable tracking~\cite{goujaudPEPitComputerassistedWorstcase2024}.
Substituting these expressions into later operations avoids auxiliary variables in the verification problem.
The MIP backend converts the remaining steps to constraints, while the simulator evaluates them on sampled inputs.
Code generation extracts the coefficient matrices and operations needed to implement the iteration.

\paragraph{Code generation for embedded optimization.}
Code generation is a key tool for producing solvers specialized to a problem family and target hardware.
Embedded applications often require static memory allocation, minimal external library dependencies, and predictable runtime.
CVXGEN~\cite{mattingleyCVXGENCodeGenerator2012}, an early code generator for small parametric quadratic programs, was used by SpaceX to generate flight code for precision rocket landing~\cite{blackmoreAutonomousPrecisionLanding2016}.
Generators for optimal control exploit the multistage structure of subproblems arising in model predictive control~\cite{ciminiEmbeddedModelPredictive2021,jerezEmbeddedOnlineOptimization2014}.
The OSQP code generator~\cite{banjacEmbeddedCodeGeneration2017} exploits fixed problem matrices to reuse a factorization across iterations of its alternating direction method of multipliers (ADMM) solver.
CVXPYgen~\cite{schallerEmbeddedCodeGeneration2022} automates the otherwise manual mapping from a CVXPY problem's parameters to the canonical data required by an embedded solver.

More general frameworks support code generation for numerical functions and configurable optimal-control algorithms.
CasADi~\cite{anderssonCasadi2019} generates C code from expression graphs, and acados~\cite{verschuerenAcados2022} combines these functions with interchangeable solver components.
\fomtools{} adds an agent that profiles and rewrites the generated solver for the target hardware, specializing its kernels to the specified iteration and problem data.
The same algorithm description also defines the performance analysis, including when the problem is supplied through CVXPYgen.

\paragraph{AI-assisted code generation.}
Recent work uses LLM coding agents to improve numerical software through execution and timing feedback~\cite{pressAlgoTuneCanLanguage2025,zimmerAgenticResearcherPractical}.
AlgoTune evaluates implementations through a Python interface and permits compiled code, for example through Cython or Numba.
Its task-specific tests validate the final outputs, allowing the agent to change the underlying algorithm.
Chen and Yin~\cite{chenTrainingNumericalIntelligence2026} use residual histories and other numerical diagnostics to guide literature search and construct reusable solver components.

\fomtools{} generates C or CUDA solvers that run without Python and compute the iterates of the specified method.
It generates the numerical data, matrix routines, and reference calculations directly from the algorithm description.
Before repeating an extracted iteration, the package compares its trajectory with the recorded algorithm, iteration by iteration, on sampled inputs.
After each code-generation phase, it checks the generated solver's outputs against the reference execution on sampled inputs.
These separate checks help distinguish errors in extracting the iteration from errors introduced during translation or code optimization.
Every generated solver passes these comparisons, although they are finite tests rather than a proof of equivalence.

\section{First-order methods and the verification problem}\label{sec:background}

This section defines the parametric optimization problem, the first-order methods, and the worst-case and probabilistic verification questions.

\paragraph{Parametric optimization.}
We consider the family of convex quadratic programs
\begin{equation}\label{prob:parametric_qp}
	\begin{array}{ll}
		\text{minimize} & (1/2) z^\tpose P z + q(x)^\tpose z\\
		\text{subject to} & Az \le b(x),
	\end{array}
\end{equation}
with variable $z \in \reals^n$ and parameter $x \in \mathcal{X} \subseteq \reals^p$, where $\mathcal{X}$ is a bounded polyhedron.
The matrices $P \in \symm_+^n$ and $A \in \reals^{m \times n}$ are fixed, and the vectors $q(x) \in \reals^n$ and $b(x) \in \reals^m$ depend affinely on $x$.
Thus the problem is convex for every $x \in \mathcal{X}$; it is a linear program when $P=0$.
Equality constraints can be written as pairs of inequalities.
We assume that a solution exists for every $x \in \mathcal{X}$.
A point $z$ is optimal if and only if there is a multiplier $y \in \reals^m$ satisfying the Karush--Kuhn--Tucker (KKT) conditions
\begin{equation}\label{eq:parametric_kkt}
    Pz + q(x) + A^\tpose y = 0, \quad Az \le b(x), \quad y \ge 0, \quad y^\tpose(Az-b(x)) = 0.
\end{equation}
For a pair $(z, y)$ with $y \ge 0$, the primal and dual residuals
\begin{equation}\label{eq:residuals}
    r_\mathrm{prim}(z, x) = \left( Az - b(x) \right)_+, \qquad r_\mathrm{dual}(z, y, x) = Pz + q(x) + A^\tpose y
\end{equation}
measure the violation of the first two conditions in~\eqref{eq:parametric_kkt}, where $(\cdot)_+$ is the elementwise positive part.

\paragraph{First-order methods.}
We solve each instance of~\eqref{prob:parametric_qp} with an iteration
\begin{equation}\label{eq:fom_iteration}
    s^{k+1} = T_k(s^k, x), \quad k = 0, 1, \ldots,
\end{equation}
where $s^k \in \reals^d$ is the algorithm state, which may contain primal, dual, and auxiliary variables.
The update $T_k : \reals^d \times \reals^p \to \reals^d$ may depend on $k$ through a prescribed step-size or momentum schedule.
For stationary methods, we write $T_k=T$.
For the methods considered here, the fixed-point equation $s^\star=T(s^\star,x)$ encodes the optimality conditions~\eqref{eq:parametric_kkt}~\cite[Secs.~2.4.3, 2.7.2]{ryuLargeScaleConvexOptimization2022}.
A primal solution can be recovered from $s^\star$.
The initial state lies in a bounded polyhedron $\mathcal{S} \subseteq \reals^d$.
A zero cold start has $\mathcal{S}=\{0\}$; for warm starts, $\mathcal{S}$ bounds the admissible initial states.
A warm start can also depend on the parameter through an affine or PWA map $s^0 = W(x)$, which the user writes as an expression of $x$.

We require each update to be a composition of affine and piecewise affine (PWA) steps in the state and parameter~\cite{ranjanVerificationFirstorderMethods2025,ranjanExactVerificationFirstOrder2025}.
Affine steps include gradient computations for quadratic objectives, momentum updates, and solves with fixed matrices.
PWA steps include box projections, soft-thresholding, rectified linear units (ReLUs), and saturated linear units.
This representation covers proximal gradient methods, accelerated variants, ADMM, and the primal-dual hybrid gradient (PDHG) method for the problem families considered here.
Since these steps are continuous and $\mathcal{S}$ and $\mathcal{X}$ are bounded, every iterate $s^k$ ranges over a bounded set, which the MIP formulation of \Sec\ref{subsec:mip_backend} needs for its big-$M$ constants.

\paragraph{Worst-case and probabilistic verification.}
Let $\phi^K(s^0,x)$ measure the performance of the algorithm after $K$ iterations, with smaller values indicating better performance.
One example is the fixed-point residual given by $\|s^K-s^{K-1}\|$.
For several methods considered below, scaling this difference by an appropriate matrix $H$ gives a residual of the KKT conditions, so we use $\|H(s^K-s^{K-1})\|$~\cite{ranjanExactVerificationFirstOrder2025}.
The scaling depends on the method and its state representation.
Another metric is the distance $\|s^K-s^\star(x)\|$ to a fixed point $s^\star(x)$.
The \emph{verification problem} (VP)~\cite{ranjanVerificationFirstorderMethods2025,ranjanExactVerificationFirstOrder2025} is
\begin{equation}\label{eq:vp}
	\begin{array}{ll}
		\text{maximize} & \phi^K(s^0, x)\\
		\text{subject to} & s^{k+1} = T_k(s^k, x), \quad k = 0, \ldots, K-1,\\
		& s^0 \in \mathcal{S}, \quad x \in \mathcal{X},
	\end{array}
\end{equation}
with variables $x$, $s^0$, and the trajectory $s^1,\ldots,s^K$.
Its optimal value is the worst-case performance after $K$ iterations.
An upper bound of $\epsilon$ on this value guarantees $\phi^K(s^0,x) \le \epsilon$ for every admissible parameter and initial state.

For probabilistic verification, we instead specify a joint distribution on $\mathcal{S}\times\mathcal{X}$ and consider the probability of exceeding a tolerance $\epsilon$,
\begin{equation}\label{eq:exceedance_probability}
    r_K(\epsilon) = \prob\bigl(\phi^K(s^0,x)>\epsilon\bigr).
\end{equation}
A bound $r_K(\epsilon)\le\alpha$, with $\alpha\in(0,1)$, means that the algorithm meets the tolerance after $K$ iterations with probability at least $1-\alpha$ under this distribution.
Independent samples of $(s^0,x)$ allow us to bound $r_K(\epsilon)$ with a chosen confidence level~\cite{sambharyaDataDrivenPerformanceGuarantees2025}.
These bounds describe the fraction of instances that exceed the tolerance, rather than the largest possible value of the metric.

The worst-case VP is nonconvex and \NPhard{} in general~\cite{ranjanVerificationFirstorderMethods2025}.
For the affine/PWA iterations and norm metrics considered here, it admits an exact mixed-integer formulation~\cite{ranjanExactVerificationFirstOrder2025}, implemented by the MIP backend (\Sec\ref{subsec:mip_backend}).
Semidefinite relaxations also bound the worst-case VP~\cite{ranjanVerificationFirstorderMethods2025}; Appendix~\ref{app:sdp_backend} describes the experimental SDP backend.
The simulator evaluates sampled trajectories and computes probabilistic bounds on~\eqref{eq:exceedance_probability} (\Sec\ref{subsec:simulator_backend}).

\subsection{Running example: Lasso, ISTA, and FISTA}\label{subsec:lasso_example}

The $\ell_1$-regularized least squares (Lasso) problem~\cite{tibshiraniRegressionShrinkageSelection1996a} is
\begin{equation}\label{prob:lasso}
	\begin{array}{ll}
		\text{minimize} & (1/2) \| D z - x \|_2^2 + \lambda \|z\|_1,
	\end{array}
\end{equation}
where $z \in \reals^n$ is the optimization variable, $D \in \reals^{m \times n}$ is a fixed dictionary, and $\lambda>0$ is a fixed regularization weight.
Introducing an auxiliary vector $t$ and replacing $\|z\|_1$ by $\ones^\tpose t$ with $-t\le z\le t$ gives an equivalent QP of the form~\eqref{prob:parametric_qp}.
The input signal $x$ ranges over a specified set $\mathcal{X}\subseteq\reals^m$ for verification.
Simulations sample $x$ from a specified distribution on $\mathcal{X}$; the Lasso experiments below use uniform draws from a box.
For a threshold $\tau\ge0$, the soft-thresholding operator acts componentwise as
\begin{equation*}
    (\softthresh_\tau(v))_i = \sign(v_i)\max\{|v_i|-\tau,0\}.
\end{equation*}
The iterative soft-thresholding algorithm (ISTA)~\cite{parikhProximalAlgorithms2014} solves~\eqref{prob:lasso} with the updates
\begin{equation}\label{eq:ista_iterations}
    s^{k+1} = \softthresh_{\lambda \eta}\!\left(s^k - \eta D^\tpose(D s^k - x)\right),
\end{equation}
for a sufficiently small step size $\eta > 0$.
This is a single PWA step with an affine argument, since we can rewrite~\eqref{eq:ista_iterations} as
\begin{equation*}
    s^{k+1} = \softthresh_{\lambda \eta}\!\left((I - \eta D^\tpose D) s^k + \eta D^\tpose x\right),
\end{equation*}
which makes the affine argument of the soft-thresholding step explicit.

The fast iterative soft-thresholding algorithm (FISTA)~\cite{beckFastIterativeShrinkageThresholding2009a} adds an auxiliary sequence $w^k$ with $w^0 = s^0$, together with a precomputed scalar sequence $\beta_{k+1} = (1 + \sqrt{1 + 4 \beta_k^2}) / 2$ starting from $\beta_0 = 1$.
The iterations are
\begin{equation}\label{eq:fista_iterations}
    \begin{array}{rl}
        s^{k+1} &= \softthresh_{\lambda \eta}\!\left(w^k - \eta D^\tpose(D w^k - x)\right),\\
        w^{k+1} &= s^{k+1} + \gamma_{k+1} (s^{k+1} - s^k),
    \end{array}
\end{equation}
where $\gamma_{k+1} = (\beta_{k} - 1) / \beta_{k+1}$.
The pair $(s^k,w^k)$ forms the algorithm state in~\eqref{eq:fom_iteration}.
For $k\ge1$, substituting the affine expression for $w^k$ into the first update gives
\begin{equation}\label{eq:fista_compact_rewrite}
    s^{k+1} = \softthresh_{\lambda \eta}\!\left((1 + \gamma_k)(I - \eta D^\tpose D) s^k - \gamma_k (I - \eta D^\tpose D) s^{k-1} + \eta D^\tpose x\right).
\end{equation}
This substitution eliminates the sequence $w^k$, and reduces both the number of variables and the number of constraints in the VP~\eqref{eq:vp}.
\fomtools{} performs these substitutions as the user records the iterations~\eqref{eq:fista_iterations}, producing~\eqref{eq:fista_compact_rewrite} without manual algebra.

\section{The \fomtools~package}\label{sec:fomtools_package}

We show how \fomtools{} represents an iteration and uses it for verification computations and code generation.
The FISTA example illustrates the connection between the Python description and the resulting VP.

\subsection{Architecture overview}\label{subsec:architecture}

\begin{figure}[tbp]
\centering
\ifpreprint
    \def\archfont{\footnotesize}\def\archhdr{\normalsize\bfseries}
    \def\archxb{5.9}\def\archxc{9.8}
    \def\archboxw{21mm}\def\archgrouptw{30mm}\def\archgroupw{34mm}
    \def\archsep{3mm}
\fi
\newlength{\archcanontw}\settowidth{\archcanontw}{\archfont Integer Program}
\begin{tikzpicture}[
    >={Latex[length=2.2mm]},
    box/.style={draw=fomtoolsblue, line width=0.9pt, rounded corners=4pt,
                fill=white, align=center, font=\archfont,
                minimum width=\archboxw, minimum height=9mm},
    canon/.style={box, text width=\archcanontw},
    group/.style={box, text width=\archgrouptw, minimum width=\archgroupw},
    grpbox/.style={rounded corners=6pt, inner sep=\archsep, line width=0.6pt},
    hdr/.style={font=\archhdr},
    arr/.style={->, draw=black!75, line width=0.7pt}
]
\node[box] (fom)   at (0,-1.0) {First-order\\Method};
\node[box] (param) at (0,-2.5) {Parameters};
\node[box] (init)  at (0,-4.0) {Initial\\Iterates};

\node[canon] (mip) at (\archxb,-1.0)  {Exact Mixed-\\Integer Program};
\node[canon] (sim) at (\archxb,-2.5)  {Simulator};
\node[canon] (cg)  at (\archxb,-5.2)  {Code\\Generation};

\node[group] (mipsolv) at (\archxc,-1.0)  {MIP solvers:\\\mbox{Gurobi}, \mbox{HiGHS}};
\node[group] (simbk)   at (\archxc,-2.5)  {Simulation:\\\mbox{NumPy}, \mbox{CuPy}};
\node[group] (cgtgt)   at (\archxc,-5.2)  {Targets:\\\mbox{Embedded C},\\\mbox{OpenMP}, \mbox{CUDA}};

\begin{scope}[on background layer]
    \node[grpbox, fill=cbBlue!12,   draw=cbBlue!55]   (userbg)  [fit=(fom)(param)(init)] {};
\node[grpbox, fill=cbOrange!16, draw=cbOrange!70]  (verifbg) [fit=(mip)(sim)(mipsolv)(simbk)] {};
    \node[grpbox, fill=cbPurple!16, draw=cbPurple!70]  (genbg)   [fit=(cg)(cgtgt)] {};
\end{scope}
\node[hdr, color=cbBlue,   above=1.5mm of userbg]  {User Input};
\node[hdr, color=cbOrange, above=1.5mm of verifbg] {Analysis};
\node[hdr, color=cbPurple, above=1.5mm of genbg]   {AI-assisted Generation};

\coordinate (vpx) at ($(userbg.east)!0.5!(verifbg.west)$);
\node[box] (vp) at (vpx |- param) {Verification\\Problem};

\draw[arr] (fom.east)   -- (vp.160);
\draw[arr] (param.east) -- (vp.180);
\draw[arr] (init.east)  -- (vp.200);

\draw[arr] (vp.east) -- (mip.west);
\draw[arr] (vp.east) -- (sim.west);
\draw[arr] (vp.east) -- (cg.west);

\draw[arr] (mip.east) -- (mipsolv.west);
\draw[arr] (sim.east) -- (simbk.west);
\draw[arr] (cg.east)  -- (cgtgt.west);
\end{tikzpicture}
\caption{One algorithm description supplies worst-case verification, sampled performance analysis, and code generation.
The code generation procedure combines deterministic data and helper generation with agent translation and optimization for the selected target.}
\label{fig:fomtools_flowchart}
\end{figure}
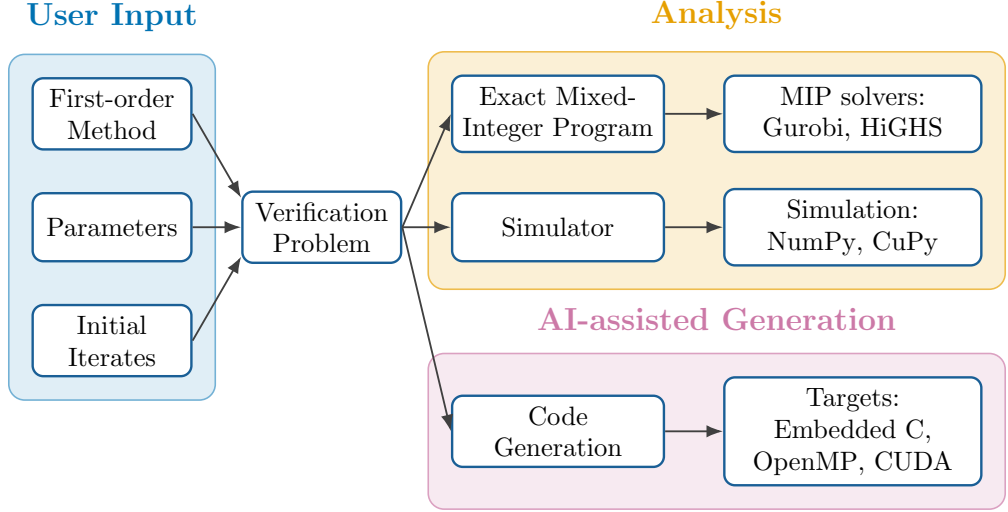

The user writes an iteration using expressions for parameters and algorithm states.
The package records the PWA steps and substitutes affine expressions to eliminate auxiliary variables.
From this representation, the MIP backend constructs constraints for worst-case verification, while the simulator evaluates sampled trajectories (Figure~\ref{fig:fomtools_flowchart}).
The same recorded iteration supplies the code generation procedure of \Sec\ref{sec:codegen}.
Appendix~\ref{app:implementation_notes} describes the interfaces between these components.

\subsection{Example: verifying FISTA}\label{subsec:workflow}

The following code bounds the $\ell_\infty$-norm of the difference between successive FISTA iterates in~\eqref{eq:fista_iterations}.
The signal ranges over $x\in[0,1]^m$, and the initial iterate is $s^0=0$.
We use $\lambda=0.1$, step size $\eta=1/\|D\|_2^2$, and a small randomly generated dictionary for illustration.

\codelisting[language=fomPython, numbers=left, firstline=1, lastline=18]{code/fista.py}

The loop records each FISTA update and solves a VP for the residual after that iteration.
Thus \texttt{results[k]} bounds the residual after $k+1$ iterations over the specified signal box.

\codelisting[language=fomPython, numbers=left, firstline=20, firstnumber=20]{code/fista.py}

The momentum expression \texttt{w} is an affine combination of \texttt{s[k+1]} and \texttt{s[k]}.
Substituting this expression into the next update gives the compact formulation~\eqref{eq:fista_compact_rewrite} without introducing a variable for \texttt{w}.

\subsection{Linear expressions}\label{subsec:linear_expressions}

\fomtools{} represents affine expressions as $Ay+b$, where $y$ stacks the \emph{leaf variables}: parameters, initial iterates, and step outputs that have not been eliminated by substitution.
Expressions support addition, scalar multiplication, and left multiplication by NumPy arrays or SciPy sparse matrices.
For example, given $A\in\reals^{n\times n}$ and $B\in\reals^{n\times m}$, the code

\codelisting[language=fomPython, numbers=left, firstline=1, lastline=5]{code/snippets.py}
represents \texttt{y} directly as an affine map of \texttt{z0} and \texttt{x}, without introducing another optimization variable.
The last line restricts $x$ to the probability simplex.
More generally, the comparison operators \texttt{==}, \texttt{<=}, and \texttt{>=} impose polyhedral constraints on expressions.
Interval arithmetic propagates their componentwise bounds through affine operations~\cite[Sec.\ 3]{gowalScalableVerifiedTraining2019}.

\subsection{Steps}\label{subsec:steps}

A \emph{step} is a primitive operator that maps an affine input to a possibly PWA output, and is the unit at which the package interfaces with each backend. 
\fomtools{} ships with four step types: affine, ReLU, soft-thresholding, and saturated linear.
The three PWA steps, ReLU, soft-thresholding, and saturated linear, each create a new leaf and record the activation pattern that the canonicalizer needs.
Explicit affine operations, such as gradient evaluations on quadratics and momentum updates, remain in the affine expressions and introduce no new leaves.
In contrast, implicit affine steps, such as linear system solves, do introduce new leaves.
The MIP backend and the simulator handle the three PWA steps natively, so the user writes the same code for every backend.

\subsection{Proxable functions}\label{subsec:proxable_functions}

A \emph{proxable function} in \fomtools{} is a convex function $f$ whose proximal operator
\begin{equation}\label{eq:prox_def}
    \prox_{\gamma f}(v) = \argmin_z \left\{ f(z) + (1/2\gamma) \|z - v\|_2^2 \right\}
\end{equation}
the package can encode explicitly, for a step size $\gamma > 0$~\cite{parikhProximalAlgorithms2014}.
The supported functions, listed in Appendix~\ref{app:prox_api}, fall into two groups.
Most have a closed form in the PWA primitives of \Sec\ref{subsec:steps}, so they compile to one step plus affine arithmetic: the $\ell_1$-norm becomes a single soft-thresholding step, the indicator of the nonnegative orthant a single ReLU, and the indicator of a box a single saturated linear step.
The rest either require a linear solve, such as convex quadratics, whose proximal operator is affine, or have no closed form in the primitives, such as the indicators of a polyhedron, a simplex, or an $\ell_1$-ball.
For these the package encodes the optimality conditions of~\eqref{eq:prox_def} directly, as linear equalities together with a complementarity condition on the multipliers.
The user applies a proxable function with \texttt{prox(f, expr, gamma)}.
For an indicator function of a set $C$, the user can also use \texttt{project(C, expr)}, since projection does not depend on $\gamma$.

\codelisting[language=fomPython, numbers=left, firstline=1, lastline=13, firstnumber=1]{code/prox.py}
Both entry points return an expression, so a proximal gradient step is a single line of user code and the decomposition into steps stays hidden.

\subsection{Verification objectives}\label{subsec:objectives}
The VP~\eqref{eq:vp} can maximize the norm of an arbitrary expression, so the user can define custom performance metrics beyond the standard fixed-point residual $\|s^{k+1} - s^k\|$. 
Examples are scaled residuals, residuals involving the problem parameters, and distances to a reference expression.
The package supports four objectives: the $\ell_1$-norm $\|v\|_1$, the squared $\ell_2$-norm $\|v\|_2^2$, the $\ell_\infty$-norm $\|v\|_\infty$, and a zero objective for feasibility checks.
The user sets the objective with \texttt{maximize}, passing \texttt{ft.norm(expr, p)} for the $\ell_1$- or $\ell_\infty$-norm, \texttt{ft.sum\_squares(expr)} for the squared $\ell_2$-norm, or \texttt{0} for a feasibility check.
The function \texttt{ft.norm} also accepts a list of expressions, which the package stacks into a single vector before applying the norm.
Stacking is useful for composite metrics such as the primal and dual residuals~\eqref{eq:residuals}, which solvers such as OSQP~\cite{stellatoOSQPOperatorSplitting2020} use together as a stopping test.
For $q(x)=q+Cx$ and $b(x)=b+Ex$, we write

\codelisting[language=fomPython, numbers=left, firstline=7, lastline=10, firstnumber=1]{code/snippets.py}

The list interface stacks both residual vectors and maximizes the $\ell_\infty$-norm of the result.
Both residuals involve the problem parameter $x$, which is permitted because the objective operates on arbitrary expressions.

\subsection{Fixed-point encoding}\label{subsec:fixed_point_encoding}

Distance-to-optimality metrics require the VP to encode an optimal solution $z^\star$ of the parametric problem.
For ISTA, the fixed-point equation $z^\star=T(z^\star,x)$ encodes the optimality conditions of the Lasso problem.
The user provides the operator as a Python function, together with bounds on $z^\star$.

\codelisting[language=fomPython, numbers=left, firstline=12, lastline=19, firstnumber=1]{code/snippets.py}

Internally, \texttt{fixed\_point} creates a new iterate $z^\star$ with the given bounds, applies the operator to obtain $T(z^\star, x)$, and adds the equality constraint $z^\star = T(z^\star, x)$.
The last line maximizes the squared distance from \texttt{zk} to a fixed point within the supplied bounds.
The bounds on $z^\star$ must hold for all $x \in \mathcal{X}$.
If the solution is not unique, the VP maximizes over all fixed points within these bounds.

\section{Analysis backends}\label{sec:backends}

The MIP backend converts the recorded iteration and performance metric into an optimization problem for worst-case verification.
We call this translation \emph{canonicalization}, following the terminology of~\cite{agrawalRewritingSystemConvex2018a,diamondCVXPYPythonEmbeddedModeling2016a}.
\fomtools{} provides the MIP backend, which constructs an exact formulation of the VP~\eqref{eq:vp} (\Sec\ref{subsec:mip_backend}), and the simulator, which evaluates the method on samples and provides statistical bounds on performance (\Sec\ref{subsec:simulator_backend}).
The experimental SDP backend relaxes the worst-case VP (Appendix~\ref{app:sdp_backend}).
All backends support the $\ell_1$, squared $\ell_2$, and $\ell_\infty$ objectives of \Sec\ref{subsec:objectives}.
The user selects one at \texttt{Verifier} construction time through the \texttt{method} keyword, and the MIP solver through the \texttt{solver} keyword, which currently accepts Gurobi and HiGHS.

\subsection{Exact verification}\label{subsec:mip_backend}

The MIP backend encodes the VP~\eqref{eq:vp} as a mixed-integer optimization problem whose optimal value equals the true worst-case performance~\cite{ranjanExactVerificationFirstOrder2025}.
The $\ell_1$ and $\ell_\infty$ objectives give MILPs, while the squared $\ell_2$ objective requires a quadratic-capable backend such as Gurobi.
The solver's best upper bound bounds the worst-case metric; a feasible incumbent gives a particular trajectory and a lower bound on the VP optimal value.
Our canonicalization backend encodes each PWA step with binary variables and big-$M$ constraints derived from the variable bounds.
For example, consider the ReLU constraint $w = \relu(v)$ for scalars $v \in \reals,~w \in \reals_+$ with bounds $l \le v \le u$.
The canonicalization introduces a binary variable $\delta \in \{0, 1\}$ and the constraints
\begin{equation*}
    w \ge 0, \quad w \ge v, \quad w \le v - l(1 - \delta), \quad w \le u \delta.
\end{equation*}
Tight variable bounds are critical for performance, since loose bounds give large big-$M$ constants, weak linear programming relaxations, and slow branch-and-bound solves.
The package also supports the convex-hull formulations of~\cite{ranjanExactVerificationFirstOrder2025}, which give tighter relaxations at the cost of additional inequality constraints.
The next three techniques produce the tight bounds that these formulations need.

\paragraph{Optimization-based bound tightening.}
The package tightens bounds through optimization-based bound tightening (OBBT)~\cite{gleixnerThreeEnhancementsOptimizationbased2017}.
For each index $i$ of a variable $v$, OBBT solves
\begin{equation*}
    \underline{v}_i = \min\; v_i \quad \text{and} \quad \overline{v}_i = \max\; v_i,
\end{equation*}
subject to the VP constraints, and updates the bounds to $[\underline{v}_i, \overline{v}_i]$.
The OBBT handler maintains a separate linear programming (LP) relaxation, encoding the continuous envelopes of the same PWA steps used by the main canonicalizer.
The package triggers OBBT automatically as new constraints and steps are added.
The tightened bounds then propagate to the expression objects for subsequent interval arithmetic, to the main canonicalizer for the big-$M$ constants, and to the OBBT handler itself for future subproblems.

\paragraph{Convex-hull cuts.}
The package also supports convex-hull cut generation~\cite[Theorem 1]{tjandraatmadjaConvexRelaxationBarrier2020}~\cite[Theorem 2.2, 2.4]{ranjanExactVerificationFirstOrder2025}.
These cuts identify violated convex-hull inequalities for the PWA step relaxations, and add them as constraints that further tighten the formulation. 
The user enables cut generation through the \texttt{convex\_hull\_cuts} keyword argument to \texttt{solve}.

\paragraph{Theory-informed bounds.}
When a theoretical convergence rate is known for the first-order method at hand, we can incorporate it into the VP~\cite{ranjanVerificationFirstorderMethods2025}.
Given a rate function $C(k)$ that bounds the residual at iteration $k$, the user calls \texttt{theory\_bound(k, target\_expr, ref\_expr)} at each iteration.
The package expands the reference expression's box by $C(k)$ in each coordinate, constrains the target to this expanded box, and propagates the tightened bounds.

An important case is an \emph{operator-informed} bound, which follows from the structure of the operator $T$ alone and therefore holds uniformly over the parameter set.
For ISTA~\eqref{eq:ista_iterations} in the running example with $\eta \le 1/\|D\|_2^2$, let $R \ge \|s^0 - z^\star(x)\|_2$ for all $s^0 \in \mathcal{S}$ and $x \in \mathcal{X}$.
Then~\cite[Theorem 2.2]{kimAnotherLookFast2018} gives, for every $k \ge 2$ and without strong convexity,
\begin{equation}\label{eq:ista_operator_bound}
    \|s^{k+1} - s^k\|_2 \le C(k) = \frac{2R}{\sqrt{(k-1)(k+2)}}.
\end{equation}
The bound holds in particular for each component, so \texttt{theory\_bound(k, s[k+1], s[k])} restricts $s^{k+1}$ to the box $[l^k - C(k)\ones,\, u^k + C(k)\ones]$, where $[l^k, u^k]$ are the current bounds on $s^k$.
The package intersects this box with the bounds already known for $s^{k+1}$, keeping the tightest bound elementwise as in~\cite[Sec.~3.1]{ranjanExactVerificationFirstOrder2025}.
The bounds on later iterates then become progressively tighter, which shrinks the big-$M$ constants and leads to faster solves with stronger bounds.

\subsection{Simulation}\label{subsec:simulator_backend}

A valid upper bound from the MIP backend bounds the worst-case metric over the specified sets, but its computational cost grows rapidly with the iteration count $K$ and the problem dimension.
The simulator backend is meant for larger problems and users who want a quick empirical assessment.
It evaluates the first-order method on a large sample of parameters and initial iterates, and reports an empirical distribution of the performance metric together with a confidence bound on the probability of exceeding a prescribed tolerance under the sampling distribution.

\paragraph{Batched evaluation.}
The simulator uses the same user-facing API as the MIP backend.
The user creates a \texttt{Verifier} with \texttt{method="simulator"}, declares parameters and initial iterates, and writes the algorithm identically.
The only difference is that the parameters and initial iterates carry user-supplied sample arrays of shape $(N, d)$, where $N$ is the number of samples.
Internally, the simulator stacks all samples into matrices over a sample dimension, and executes each iteration as a batched operation in single precision.
Affine maps become matrix multiplications, and PWA steps become element-wise operations, \eg, \texttt{maximum(result, 0)} for ReLU.
The simulator evaluates the objective per sample at each iteration, which yields a residual matrix of shape $(N, K)$.
After the solve, the user can query summary statistics through a dictionary interface: the mean, standard deviation, minimum, maximum, median, and quantiles of the residual, together with the empirical risk, the risk bound, the runtime, and the throughput.

\paragraph{GPU execution.}
When a CUDA-compatible GPU is available, the simulator uses CuPy~\cite{okutaCuPyNumPyCompatibleLibrary2017} as a drop-in replacement for NumPy, with sparse coefficient matrices stored in compressed sparse row (CSR) format.
Otherwise, it falls back to NumPy and SciPy on the central processing unit (CPU), with no change to the user code.
The GPU path is effective because the workload is highly parallelizable.
Each sample is independent, and the per-sample operations are identical~\cite{blinBatchedFirstOrderMethods2026}.
For accurate timing, the package performs a warm-up execution to trigger lazy kernel compilation, synchronizes the CUDA stream, and then times the simulation, reporting throughput in samples per second alongside the statistical summaries.

\paragraph{Statistical guarantees.}
Fix the algorithm, iteration count $K$, tolerance $\epsilon \ge 0$, and sampling distribution from \Sec\ref{sec:background} before drawing the evaluation sample.
For $N$ independent pairs $(s_i^0,x_i)$ from this distribution, let $\phi_i^K=\phi^K(s_i^0,x_i)$.
The simulator estimates the exceedance probability~\eqref{eq:exceedance_probability} by the empirical fraction
\begin{equation*}
    \hat r_K(\epsilon) = \frac{1}{N}\sum_{i=1}^N \mathbf{1}\{\phi_i^K>\epsilon\}.
\end{equation*}
With probability at least $1-\delta$ over the sample, the Kullback--Leibler (KL) inverse bound~\cite{langfordNotBoundingTrue2001,sambharyaDataDrivenPerformanceGuarantees2025} gives
\begin{equation*}
    r_K(\epsilon) \le \mathrm{kl}^{-1}\!\left(\hat r_K(\epsilon) \;\middle|\; \frac{\log(2/\delta)}{N}\right),
\end{equation*}
where $\delta\in(0,1)$ is the confidence failure probability.
Here $\mathrm{kl}^{-1}(q \mid c)$ is the largest $p \in [q,1]$ satisfying $q\log(q/p)+(1-q)\log((1-q)/(1-p))\le c$, computed numerically.
The default is $\delta=0.05$; both $\delta$ and $\epsilon$ can be chosen by the user.

To obtain bounds at several iteration counts or tolerances, use $\delta/M$ in the formula above, where $M$ is the total number of bounds computed.
A union bound then makes all of them hold with probability at least $1-\delta$.
The algorithm choices, iteration counts, and candidate tolerances must be fixed independently of the evaluation sample.
For a joint guarantee across several algorithms, include their bounds in $M$.

For a fixed tolerance, we can plot the risk bound against the iteration count.
We can also seek a tolerance satisfied by at least a fraction $1-\alpha$ of instances, for a chosen $\alpha \in (0,1)$.
At each iteration, select the smallest candidate tolerance whose risk bound is at most $\alpha$, if one exists.
With probability at least $1-\delta$, the selected tolerances meet this requirement at all computed iterations~\cite[Appendix~B.2]{sambharyaDataDrivenPerformanceGuarantees2025}.
For example, $\alpha = 0.05$ gives these probabilistic upper bounds on the $95$th percentile of the residual.
The candidate tolerances can be chosen using an independent reference sample.

\section{Code generation}\label{sec:codegen}
\fomtools{} generates a standalone solver from the same iteration used for verification and simulation.
It first prepares the fixed numerical data and supporting routines.
A coding agent then translates the iteration into code and optionally optimizes the implementation.
Numerical checks compare its outputs with a NumPy reference on sampled inputs.
The result is a C or CUDA library with a Python interface.

\subsection{The solver loop}\label{subsec:codegen_ir}

The generated solver maintains the state needed by the next update: one iterate for ISTA and two for FISTA.
The problem parameters remain fixed during a solve.
\fomtools{} represents each iteration by a computation graph.
We number its steps so that each step uses only the state, the parameters, and outputs of earlier steps.
With $s$ the state, $x$ the parameters, and $t_1, \ldots, t_J$ the intermediate outputs, iteration $k$ has the form
\begin{equation}\label{eq:codegen_ir}
	\begin{array}{rl}
		t_j &= g_j\!\left(a_j^k(s, x, t_1, \ldots, t_{j-1})\right), \quad j = 1, \ldots, J,\\
		s^+ &= a_{J+1}^k(s, x, t_1, \ldots, t_J),
	\end{array}
\end{equation}
where each $a_j^k$ is affine and each $g_j$ is a PWA step (\Sec\ref{subsec:steps}), the identity, or the solution of a linear system with a fixed matrix.
An affine map need not depend on all its arguments.

For ISTA~\eqref{eq:ista_iterations}, the state $s$ is the current iterate, $x$ is the input signal, and $J=1$.
The maps are
\begin{equation*}
    \begin{aligned}
        a_1^k(s,x) &= (I-\eta D^\tpose D)s+\eta D^\tpose x,\\
        g_1(v) &= \softthresh_{\lambda\eta}(v),\\
        a_2^k(s,x,t_1) &= t_1.
    \end{aligned}
\end{equation*}
The first affine map computes a gradient step, $g_1$ applies soft-thresholding, and the final map sets the next state to $t_1$.

A \emph{block} contains the operations and coefficient data for one iteration.
\fomtools{} extracts these blocks from the recorded computation graph.
For a stationary method such as ISTA, the affine maps do not depend on $k$, so the solver repeats one block in a loop.
In code, the user records this iteration once and calls \texttt{codegen\_rolled}.

\codelisting[language=fomPython, numbers=left]{code/codegen_rolled.py}
The generated solver runs this loop up to an iteration cap set at run time.

In FISTA~\eqref{eq:fista_iterations}, the operations are unchanged but the momentum $\gamma_k$ varies with $k$.
The user supplies these values as a \emph{schedule} and records one iteration as a function of $\gamma$.
Its coefficient matrices have the form $M(\gamma) = M_0 + \gamma M_1$, with $M_0$ and $M_1$ fixed.
The solver stores these matrices once and reads $\gamma_k$ from a table at iteration $k$.
Thus FISTA uses one iteration block even though its coefficients change.
For several schedules, each coefficient matrix must depend affinely on the scheduled scalars.

A recorded run stores every iteration separately.
\fomtools{} checks whether repeating the first iteration reproduces the recorded iterates on sampled parameter values.
If it does, within a numerical tolerance, the solver is a loop that runs for the recorded number of iterations.
Otherwise, the package generates \emph{unrolled} code, which writes out each recorded iteration in sequence.
The network experiment in \Sec\ref{subsec:network} compares a loop with scheduled coefficients and unrolled code for the same method.
Appendix~\ref{app:codegen_inputs} gives the input interfaces and code examples.

\subsection{Generating and checking solver code}\label{subsec:codegen_verify}

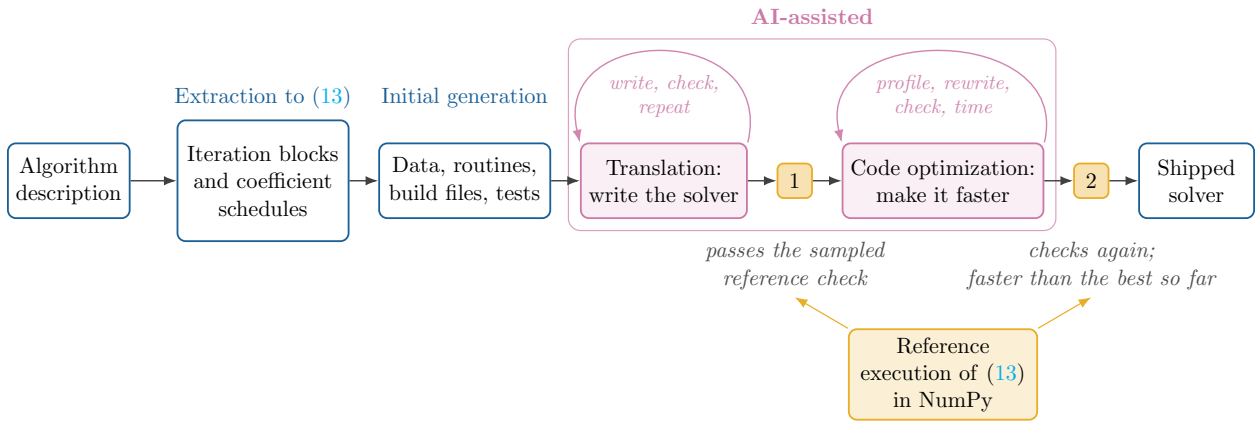
\begin{figure}[ht]
\centering
\ifpreprint
    \def\cgfont{\small}\def\cglbl{\small}\def\cgloop{\footnotesize}\def\cghdr{\small\bfseries}
\fi
\resizebox{\textwidth}{!}{\begin{tikzpicture}[
    >={Latex[length=2.2mm]},
    box/.style={draw=fomtoolsblue, line width=0.9pt, rounded corners=4pt,
                fill=white, align=center, font=\cgfont,
                minimum width=20mm, minimum height=13mm},
    agent/.style={box, draw=cbPurple, fill=cbPurple!12},
    gate/.style={draw=cbOrange!85, line width=0.9pt, rounded corners=3pt,
                 fill=cbOrange!30, align=center, font=\cgfont,
                 minimum width=6mm, minimum height=6mm, inner sep=1.5pt},
    ref/.style={draw=cbOrange!85, line width=0.9pt, rounded corners=4pt,
                fill=cbOrange!16, align=center, font=\cgfont, text width=30mm,
                minimum height=10mm},
    arr/.style={->, draw=black!75, line width=0.7pt},
    thin arr/.style={->, draw=cbOrange!85, line width=0.6pt},
    lbl/.style={font=\cglbl\itshape, align=center, text=black!70},
    aibox/.style={draw=cbPurple!70, line width=0.6pt, rounded corners=6pt, inner sep=2mm}
]
\node[box]                     (vp)   {Algorithm\\description};
\node[box, right=8mm of vp, text width=27mm, minimum height=21mm] (ir)
     {Iteration blocks\\and coefficient\\schedules};
\node[font=\cglbl, text=fomtoolsblue, above=0.5mm of ir] (exlbl) {Extraction to~\eqref{eq:codegen_ir}};
\node[box,   right=5mm of ir, text width=27mm] (supp)
     {Data, routines,\\build files, tests};
\node[font=\cglbl, text=fomtoolsblue, anchor=base] at (supp |- exlbl.base) {Initial generation};
\node[agent, right=5mm of supp] (pa)  {Translation:\\write the solver};
\node[gate,  right=5mm of pa]   (g1)  {1};
\node[agent, right=5mm of g1]   (pb)  {Code optimization:\\make it faster};
\node[gate,  right=5mm of pb]   (g2)  {2};
\node[box,   right=5mm of g2]   (out) {Shipped\\solver};

\draw[arr] (vp)   -- (ir);
\draw[arr] (ir)   -- (supp);
\draw[arr] (supp) -- (pa);
\draw[arr] (pa)   -- (g1);
\draw[arr] (g1)   -- (pb);
\draw[arr] (pb)   -- (g2);
\draw[arr] (g2)   -- (out);

\draw[arr, draw=cbPurple!90] (pa.north east) .. controls +(5mm,21mm) and +(-5mm,21mm) .. (pa.north west);
\node[lbl, above=3mm of pa, text=cbPurple!90, font=\cgloop\itshape] (lpa)
     {write, check,\\repeat};
\draw[arr, draw=cbPurple!90] (pb.north east) .. controls +(5mm,21mm) and +(-5mm,21mm) .. (pb.north west);
\node[lbl, above=3mm of pb, text=cbPurple!90, font=\cgloop\itshape] (lpb)
     {profile, rewrite,\\check, time};

\begin{scope}[on background layer]
    \coordinate (apex) at ([yshift=16mm]pa.north);  \node[aibox] (aibg) [fit=(pa)(pb)(apex)] {};
\end{scope}
\node[font=\cghdr, text=cbPurple, above=1mm of aibg] {AI-assisted};

\node[lbl, below=6mm of g1] (l1) {passes the sampled\\reference check};
\node[lbl, below=6mm of g2] (l2) {checks again;\\faster than the best so far};

\node[ref, below=19mm of pb] (pyref) {Reference execution of~\eqref{eq:codegen_ir} in NumPy};
\draw[thin arr] (pyref.north west) -- (l1.south);
\draw[thin arr] (pyref.north east) -- (l2.south);
\end{tikzpicture}}
\caption{The code generation pipeline.
Extraction obtains the iteration blocks~\eqref{eq:codegen_ir} and any coefficient schedules from the algorithm description.
The agent writes the solver loop during translation and may also rewrite supporting routines during code optimization.
The numbered gates check the generated code against fixed references.
The code optimization phase runs in rounds, and gate~2 checks each round and, with profiling enabled, times it against the fastest solver so far.
If gate~2 rejects a round, the best earlier solver is retained and the phase ends.
If translation does not pass gate~1 within its attempt and time limits, generation fails.
For CUDA without GPU access, both gates check only compilation.}
\label{fig:codegen_pipeline}
\end{figure}

\fomtools{} generates the numerical data, supporting routines, and tests from the iteration~\eqref{eq:codegen_ir}.
The coding agent first writes the solver loop using these routines.
During code optimization, it may also rewrite the matrix-vector products and linear solves, provided these routines can be checked on the target hardware.
The problem data and reference checks remain fixed.
Figure~\ref{fig:codegen_pipeline} shows the pipeline.

\paragraph{Targets and interface.}
The default target is embedded C in the style of the OSQP code generator~\cite{banjacEmbeddedCodeGeneration2017}.
The package also supports OpenMP for multicore CPUs and CUDA for GPUs.
All targets expose the same five-function interface.
Figure~\ref{fig:codegen_files} shows the files of a generated solver, and Appendix~\ref{app:codegen_files} describes the targets and the interface.

\begin{figure}[ht]
    \centering
\newcommand{\treerow}[4]{\ifnum#2>0\relax
            \draw[gray!70] ({#2*0.8-0.37}, {-#1*0.46+0.08}) -- ({#2*0.8+0.06}, {-#1*0.46+0.08});
        \fi
        \node[anchor=base west, font=\small\ttfamily] at ({#2*0.8}, {-#1*0.46}) {#3};
        \node[anchor=base west, font=\small] at (5.2, {-#1*0.46}) {#4};}
    \newcommand{\treestem}[3]{\draw[gray!70] ({#3*0.8+0.43}, {-#1*0.46-0.12}) -- ({#3*0.8+0.43}, {-#2*0.46+0.08});}
    \begin{tikzpicture}
        \fill[fomtoolsblue!10, rounded corners=2pt] (1.52, -2.42) rectangle (11.55, -1.08);
        \treestem{0}{22}{0}
        \treestem{2}{8}{1}
        \treestem{9}{14}{1}
        \treestem{15}{16}{1}
        \treestem{17}{18}{1}
        \treerow{0}{0}{solver/}{}
        \treerow{1}{1}{include/}{dimensions, workspace struct, declarations}
        \treerow{2}{1}{src/}{}
        \treerow{3}{2}{\textbf{fom\_solver.c}}{\textbf{the iteration loop}}
        \treerow{4}{2}{fom\_ops.c}{matrix-vector and per-block routines}
        \treerow{5}{2}{fom\_linsolve.c}{linear solves and their factors}
        \treerow{6}{2}{fom\_data.c}{matrix data, thresholds, schedules}
        \treerow{7}{2}{fom\_objective.c}{verification metric for stopping}
        \treerow{8}{2}{fom\_profile.c}{region profiler, compiled only for profiling}
        \treerow{9}{1}{test/}{}
        \treerow{10}{2}{verify\_main.c}{the numerical acceptance check}
        \treerow{11}{2}{verify\_data.h}{its samples, reference output, and probes}
        \treerow{12}{2}{ops\_ref.\char123 c,h\char125}{reference matrix-vector products}
        \treerow{13}{2}{kernel\_check.py}{rules for the rewritten sources}
        \treerow{14}{2}{bench\_main.c}{the timing benchmark}
        \treerow{15}{1}{python/}{}
        \treerow{16}{2}{fom\_solver.py}{wrapper over the shared library}
        \treerow{17}{1}{release/}{}
        \treerow{18}{2}{fom\_*.c}{production copies of the highlighted files}
        \treerow{19}{1}{SPEC.md}{the specification, in prose}
        \treerow{20}{1}{spec.json, data.npz}{the same specification, and the data}
        \treerow{21}{1}{OPTIMIZATIONS.md}{timings, diffs, and kernel checks}
        \treerow{22}{1}{Makefile, README.md}{build rules and usage}
    \end{tikzpicture}
    \caption{Files for a generated CPU solver.
    During translation, the agent writes only \texttt{src/fom\_solver.c}, in bold.
    During code optimization, it may rewrite all three highlighted files.
    The package restores the other generated source files, headers, build files, and tests before each check.
    Matrix-vector products and direct linear solves have additional tests against fixed references; iterative linear solves are checked through the complete solver.}
    \label{fig:codegen_files}
\end{figure}

\paragraph{The agent.}
The agent works in two phases, shown in purple in Figure~\ref{fig:codegen_pipeline}.
In the translation phase, it writes the solver loop from the iteration blocks, calling the generated matrix-vector and linear-solve routines.
It revises the code until the numerical check described below passes, or until it reaches its limit on attempts or time.
The optional code optimization phase starts from this accepted solver.
The agent profiles the code and rewrites its most expensive operations.
The phase runs in rounds, and each round starts a new agent session from the fastest verified solver so far.
The package keeps a round only if its code passes the checks and runs faster than that solver.
The phase stops after a round that fails, is not faster, or gains less than $2\%$, and after at most five rounds.
The package drives the agent through a common interface and does not depend on a particular model.
It includes drivers for Claude Code, Codex, and pi, and any other command-line coding agent can be added through a short configuration entry.

\paragraph{The numerical acceptance check.}
After each phase, the package runs the generated solver on sampled parameters and compares its final state with a \emph{reference execution} of the iteration blocks~\eqref{eq:codegen_ir} in NumPy.
The parameters are drawn uniformly from the box defined by the bounds on $x$, so the samples need not satisfy constraints that couple the coordinates, such as a budget constraint.
The check passes if the maximum absolute error over all samples and state entries is at most $\epsilon_\mathrm{gen}$.
The reference uses the extracted blocks, so for a recorded run the package also compares the extracted loop or unrolled blocks with a NumPy execution of the recorded computation graph.
When the agent may rewrite the matrix-vector products and direct linear solves, these routines are also tested individually against reference results computed at generation time.
Appendix~\ref{app:codegen_files} gives the number of samples, the tolerances, and the routine tests.

Before each check, the package restores the fixed data and test files, so the agent cannot change what its code is tested against.
After each round of code optimization, the package checks the new solver again.
With profiling enabled, it also times the new solver against the fastest solver so far, back to back.
The first round may be at most $2\%$ slower than the translated solver, and each later round must be faster than the best earlier one.
If a round fails either test, the package keeps the best solver so far, which is the translated solver when no round passed.
Without profiling, the package keeps a rewrite that passes the checks and runs a single round.
For CUDA, the numerical checks require a GPU.
Without one, the package only compiles the code and lets the agent edit only the solver loop.
These checks are finite tests.
They do not prove that the generated code agrees with the iteration for every parameter, and a VP bound does not automatically apply to the floating-point execution.

\subsection{Modeling in CVXPY}\label{subsec:codegen_cvxpygen}

Users can also generate a solver directly from a parametric CVXPY model.
CVXPYgen~\cite{schallerEmbeddedCodeGeneration2022} converts the model to a quadratic program and generates code to update its data and recover the model's variables.
\fomtools{} supplies the first-order method and the solver code.
The quadratic and constraint matrices must be fixed.
Parameters may enter only the linear objective and the constraint bounds.
The portfolio example in \Sec\ref{subsec:portfolio} uses this interface for problems with up to $3000$ assets.

The user selects a method through the \texttt{fom} argument.
The built-in choices are \texttt{admm}, the ADMM iteration of OSQP~\cite{stellatoOSQPOperatorSplitting2020} with a constant over-relaxation, given in Appendix~\ref{app:experiments_portfolio}, and \texttt{pdhg}, which applies PDHG to problems with no quadratic term.
Each choice calls a \emph{builder} that records the method on a verifier.
A user can supply another builder through the same argument, using the input forms in Appendix~\ref{app:codegen_inputs}.
The built-in methods start each solve from the zero state.

Let $\theta$ be the stacked vector of CVXPY parameters.
The solver's parameter vector $x$ stacks the entries of the linear cost and the constraint bounds that depend on $\theta$.
The generated code computes it as $x=C\theta+d$, with fixed $C$ and $d$.
This map also includes a fixed diagonal scaling of the QP, applied by default to improve its numerical conditioning.
Through \texttt{param\_bounds}, the user supplies a box $\Theta$ of possible values of $\theta$.
The package then sets
\begin{equation}\label{eq:codegen_interval_hull}
    \mathcal{X} = \{\, x \mid \underline{x} \le x \le \overline{x} \,\}, \qquad
    \underline{x}_i = \min_{\theta \in \Theta} (C\theta+d)_i, \qquad
    \overline{x}_i = \max_{\theta \in \Theta} (C\theta+d)_i.
\end{equation}
These extrema have closed-form expressions because the map is affine and $\Theta$ is a box.
The box $\mathcal{X}$ contains every parameter vector obtained from $\Theta$, but can also contain vectors outside that image, for example when one CVXPY parameter enters several entries of $x$.
The acceptance check samples from this box, so it tests the solver on a superset of the values it will see.

The interface returns the solution in the model's original units.
Optional early stopping can use primal feasibility and stationarity residuals measured in those units.
The following example generates an ADMM solver for a parametric quadratic program in which the linear cost $q$ is a parameter.
The \texttt{config} argument of the same call selects the target in \Sec\ref{subsec:codegen_verify}, and \texttt{code\_optimization=True} adds the code optimization phase.

\codelisting[language=fomPython]{code/cvxpygen.py}

\section{Numerical experiments}\label{sec:experiments}

We evaluate \fomtools{} on three parametric problem families from sparse coding, network optimization, and optimal control.
For each family, we verify the methods and generate their solvers.
A final portfolio rebalancing example compares solvers generated through the CVXPYgen front end with the OSQP solvers that CVXPYgen generates.
\Sec\ref{subsec:experiments_setup} describes the setup shared by all experiments, \Sec\ref{subsec:lasso}--\Sec\ref{subsec:opt_control} cover the three families, and \Sec\ref{subsec:portfolio} covers the portfolio example.

\subsection{Setup}\label{subsec:experiments_setup}

Appendix~\ref{app:experiments} gives the formulations, the iterations, and the problem data for all experiments.

\paragraph{Performance metric.}
Unless stated otherwise, the VP maximizes the $\ell_\infty$-norm of the scaled fixed-point residual $\|H(s^K - s^{K-1})\|_\infty$ of \Sec\ref{sec:background}.
The matrix $H$ maps the fixed-point residual to a residual of the optimality conditions of the underlying problem, so the metric has the same meaning for different methods.
Appendix~\ref{app:experiments} gives $H$ for each method, following the derivations of~\cite{ranjanExactVerificationFirstOrder2025}.

\paragraph{Probabilistic bounds.}
For each method, we draw independent reference and evaluation samples of $10^3$ parameter values uniformly from $\mathcal{X}$.
At each $K$, we use $100$ geometrically spaced tolerances from $0.9$ times the reference sample's $95$th-percentile residual to twice its largest residual.
We apply the construction of \Sec\ref{subsec:simulator_backend} with $\alpha = \delta = 0.05$, taking the union bound over all tolerances and plotted iteration counts in each curve.
For each curve, with probability at least $0.95$ over the evaluation sample, the bounds on the $95$th percentile hold at every plotted $K$.
In contrast, the MIP bound holds for every parameter in $\mathcal{X}$.

\paragraph{Verification solves.}
We solve the MILPs with Gurobi 13.0~\cite{gurobi} on $16$ CPU cores with $32$\,GB of RAM on the Della cluster at Princeton University.
The duality gap tolerance is $5\%$, and the time limit is two hours for each $K$.
Every solve in this section reaches this gap tolerance.
Building and solving the VPs up to the plotted iteration counts takes at most $1.3$ hours for each network flow variant, $0.8$ and $5.1$ hours for the restart periods $K_\mathrm{r} = 6$ and $12$ of the control family, and $1.3$ and $5.6$ hours for ISTA and FISTA on the Lasso family.
The largest models have about $4400$ variables, of which up to $2100$ are binary.
We solve the VPs for $K = 1, 2, \ldots$ in sequence, and each solve starts from the previous incumbent.
The figures show Gurobi's best bound, which is an upper bound on the VP optimal value up to the solver's numerical tolerances.

\paragraph{Code generation and timing.}
For each verified method, we also generate a C solver for the verified instance with the translation and code optimization phases of \Sec\ref{subsec:codegen_verify}, whose numerical acceptance check draws parameters from the verified parameter set.
Every code generation experiment in this paper uses Claude Opus 5.5 as the coding agent, run through the Claude Agent SDK.
Any other supported agent can be used without changing the pipeline.
Agent runs are not deterministic, so another run or another agent may find different optimizations and speedups, while every generated solver must still pass the numerical acceptance check.
We therefore generate every solver three times, in independent agent runs, and report the median over the three solvers together with the range of the speedups.
We time each solver for a complete $200$-iteration run against a direct NumPy implementation of the same iteration, in double precision.
All timed C and CUDA code, including the solvers that CVXPYgen generates, is compiled with \texttt{-O3}.
The code follows embedded C conventions, but we measure the timings on a desktop processor rather than on embedded hardware.
The fixed iteration count compares implementations of each method, not the time to reach a common accuracy.
Each family closes with a larger code generation experiment, which runs on the same machine unless stated otherwise.
Appendix~\ref{app:experiments} gives the timing protocol, and \Sec\ref{sec:introduction} gives the repositories of the package and of the experiments.

\newcommand{\codegentable}[1]{\small
    \setlength{\tabcolsep}{4pt}
    \begin{tabular}{@{}lrrrrrrr@{}}
        \toprule
        & \multicolumn{3}{c}{Time per solve ($\mu$s)} & \multicolumn{4}{c}{Speedup}\\
        \cmidrule(lr){2-4}\cmidrule(l){5-8}
        Method & NumPy & transl. & optim. & transl. & optim. & total & range\\
        \midrule
        \csvreader[head to column names, late after line=\\, filter strcmp={\group}{#1}]{./data/codegen_timing.csv}{}{\method & \numpyus & \translatedus & \optimizedus & $\speedtrans\times$ & $\speedopt\times$ & $\speedtotal\times$ & $\speedtotalmin$--$\speedtotalmax\times$}
        \bottomrule
    \end{tabular}}

\subsection{Lasso}\label{subsec:lasso}

We return to the running example of \Sec\ref{subsec:lasso_example}.
The dictionary is underdetermined, so $D^\tpose D$ is singular and the problem is not strongly convex.
The parameter set $\mathcal{X}$ is the smallest box that contains the measurements of $100$ sparse signals through this dictionary, similarly to~\cite{chenTheoreticalLinearConvergence2018}.
We compare ISTA~\eqref{eq:ista_iterations} with FISTA~\eqref{eq:fista_iterations}.
Appendix~\ref{app:experiments_lasso} gives the step sizes and the problem data.

\begin{figure}[ht]
\centering
\begin{tikzpicture}
	\begin{semilogyaxis}[
		width=0.6\textwidth,
		height=0.42\textwidth,
		xlabel={$K$},
		ylabel={KKT residual},
		xmin=1, xmax=20,
		xtick={2,4,...,20},
		ytick={0.02,0.05,0.1,0.2,0.5,1,2},
		yticklabels={$0.02$,$0.05$,$0.1$,$0.2$,$0.5$,$1$,$2$},
		yminorticks=false,
		grid=major,
		legend cell align=left,
		legend style={font=\small, at={(1.03,1)}, anchor=north west},
	]
		\addplot[color=fomtoolsblue, thick, mark=*, mark size=1.2pt] table [x=K, y=ista_mip, col sep=comma] {./data/lasso_residual.csv};
		\addlegendentry{ISTA (worst-case)}
		\addplot[color=fomtoolsblue, dashed] table [x=K, y=ista_prob95, col sep=comma] {./data/lasso_residual.csv};
		\addlegendentry{ISTA (probabilistic)}
		\addplot[color=red!70!black, thick, mark=*, mark size=1.2pt] table [x=K, y=fista_mip, col sep=comma] {./data/lasso_residual.csv};
		\addlegendentry{FISTA (worst-case)}
		\addplot[color=red!70!black, dashed] table [x=K, y=fista_prob95, col sep=comma] {./data/lasso_residual.csv};
		\addlegendentry{FISTA (probabilistic)}
	\end{semilogyaxis}
\end{tikzpicture}
\caption{Worst-case KKT residual~\eqref{eq:lasso_residual} of ISTA and FISTA on the nonstrongly convex Lasso family after $K$ iterations, computed by the MIP backend (solid), and probabilistic bounds on the $95$th percentile of the residual over sampled measurements (dashed), described in \Sec\ref{subsec:experiments_setup}.
The two methods coincide for $K \le 2$, since the first momentum coefficient is zero.}
\label{fig:lasso_residual}
\end{figure}
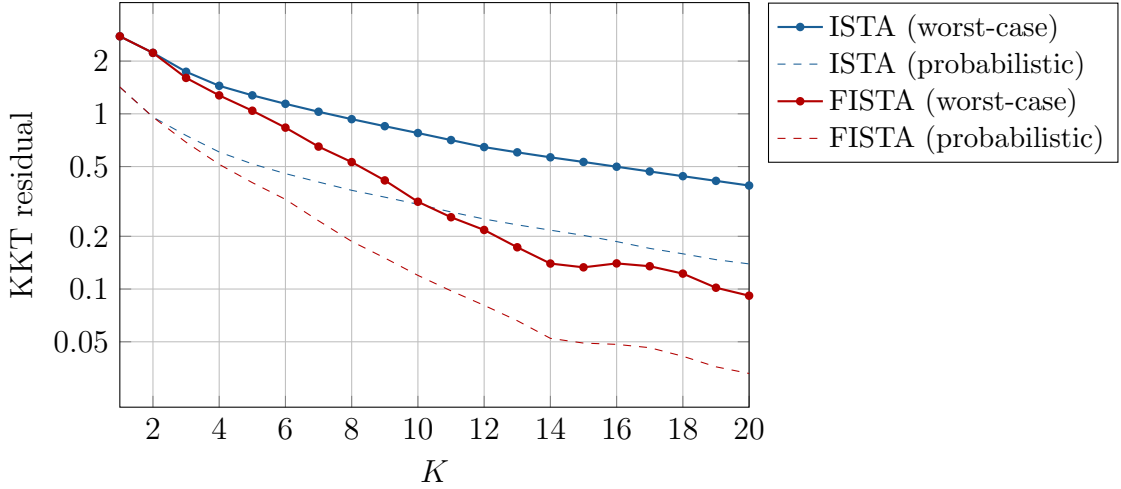

\paragraph{Results.}
Figure~\ref{fig:lasso_residual} shows that acceleration pays off in the worst case.
From $K = 3$ on, the FISTA bound lies below the best trajectory that Gurobi finds for ISTA, so FISTA provably has the smaller worst-case residual, about four times smaller at $K = 20$, where the FISTA bound is $0.092$ and Gurobi finds an ISTA trajectory with residual $0.376$.
The VP also reveals what sampling misses.
The worst case of FISTA is not monotone in $K$.
For both methods, the worst case is $1.9$ to $3.0$ times the probabilistic bound on the $95$th percentile, so the hardest measurements are rare under uniform sampling.

Table~\ref{tab:codegen_lasso} shows that code generation speeds up both methods by similar overall factors, $32$ and $37$ times.
The FISTA momentum coefficient is declared as a schedule, as in \Sec\ref{subsec:codegen_ir}, and the translated FISTA solver is more than three times slower than the translated ISTA solver.
The code optimization phase removes almost all of this overhead.

\begin{table}[ht]
    \centering
    \caption{Generated C solvers for the verified Lasso instance: time per $200$-iteration solve, and speedups of the translated solver over NumPy, of the optimized solver over the translated one, and in total.
    Each entry is the median over three independent agent runs, and the last column gives the range of the total speedup.}
    \label{tab:codegen_lasso}
    \codegentable{Lasso}
\end{table}

\paragraph{GPU code generation.}
The MIP limits the VP to small instances, but code generation scales far beyond them and also targets GPUs.
We generate ISTA and FISTA solvers for a sparse dictionary $D$ with $4 \cdot 10^4$ columns and about $260$ nonzeros per row, on the OpenMP and CUDA targets of \Sec\ref{subsec:codegen_verify}.
Appendix~\ref{app:experiments_lasso} gives the instance.
The extracted iteration multiplies by $D$ and $D^\tpose$ separately, as a hand-written solver would.
The precomputed matrix $I - \eta D^\tpose D$ is $34\%$ dense on this instance, so using it would read about $100$ times as much matrix data per iteration.
Table~\ref{tab:codegen_ista_fista} shows that the generated CUDA solvers are $4.6$ and $3.3$ times faster than a hand-written implementation in CuPy, whose sparse products call cuSPARSE.
On the CPU, the generated OpenMP solvers on $16$ threads are $29$ and $26$ times faster than a hand-written SciPy implementation on one core.
We use FISTA as an example of a method that a user has chosen to deploy, not as the best algorithm for this problem.
The agent exploits the structure of both the method and the target.
In every FISTA run on the CPU, it computes the momentum terms with a single product with $D$ per iteration, instead of one product for each term.
On the GPU, it changes how the sparse products read the iterate, for example by staging slices of the iterate in the shared memory of each block.
Table~\ref{tab:codegen_lasso_agent_time} reports the agent time: translation took $2$ to $11$ minutes per solver, and code optimization $1.7$ to $6.2$ hours, in four or five rounds.

\begin{table}[ht]
    \centering
    \caption{Generated and hand-written solvers for ISTA and FISTA with $D \in \reals^{10{,}000 \times 40{,}000}$ and about $260$ nonzeros per row, at $200$ iterations per solve in double precision.
    On the CPU, the generated OpenMP solvers use $16$ threads and the hand-written SciPy implementation one core.
    On the GPU, the generated CUDA solvers and the hand-written CuPy implementation, whose sparse products call cuSPARSE, run on the same device.
    Each speedup is that of the generated solver over the hand-written one, and the range is over three independent agent runs.
    All generated solvers passed the numerical acceptance check with tolerance $\epsilon_\mathrm{gen} = 10^{-5}$.
    Times are medians over seven interleaved rounds on an AMD EPYC 7413 CPU and an NVIDIA A100 GPU on the same node, and the generated times are medians over the three runs.
    The instance is given in Appendix~\ref{app:experiments_lasso}.}
    \label{tab:codegen_ista_fista}
    \small
    \begin{adjustbox}{max width=\linewidth}
    \begin{tabular}{@{}lrrrrrrrr@{}}
        \toprule
        & \multicolumn{4}{c}{CPU} & \multicolumn{4}{c}{GPU} \\
        \cmidrule(lr){2-5}\cmidrule(l){6-9}
        & Hand-written & Generated & & & Hand-written & Generated & & \\
        Method & SciPy (s) & OpenMP (s) & Speedup & range & CuPy (s) & CUDA (s) & Speedup & range\\
        \midrule
        \begin{filecontents*}{./data/codegen_lasso_rect.csv}
method,scipys,omps,speedomp,speedompmin,speedompmax,cupys,cudas,speedcuda,speedcudamin,speedcudamax
ISTA,0.819,0.0280,29,19,30,0.0407,0.0088,4.6,4.3,5.5
FISTA,0.820,0.0311,26,24,29,0.0471,0.0142,3.3,3.2,5.3
\end{filecontents*}
\csvreader[head to column names, late after line=\\]{./data/codegen_lasso_rect.csv}{}{\method & \scipys & \omps & $\speedomp\times$ & $\speedompmin$--$\speedompmax\times$ & \cupys & \cudas & $\speedcuda\times$ & $\speedcudamin$--$\speedcudamax\times$}
        \bottomrule
    \end{tabular}
    \end{adjustbox}
\end{table}

\begin{table}[ht]
    \centering
    \caption{Agent wall-clock time to generate the solvers of Table~\ref{tab:codegen_ista_fista}, and the number of code optimization rounds, as the median and the range over three independent agent runs.
    The phase times include the package's checks after each phase and round, and the builds and timings of the agent on GPU nodes of the cluster.
    Everything else takes at most a few seconds per solver.}
    \label{tab:codegen_lasso_agent_time}
    \small
    \begin{tabular}{@{}llrrrrrr@{}}
        \toprule
        & & \multicolumn{2}{c}{Translation (min)} & \multicolumn{2}{c}{Code optimization (h)} & \multicolumn{2}{c}{Rounds}\\
        \cmidrule(lr){3-4}\cmidrule(lr){5-6}\cmidrule(l){7-8}
        Method & Target & median & range & median & range & median & range\\
        \midrule
        \begin{filecontents*}{./data/codegen_lasso_agent_time.csv}
method,target,tr,trmin,trmax,opt,optmin,optmax,rounds,roundsmin,roundsmax
ISTA,OpenMP,1.7,1.6,1.9,3.3,1.7,3.3,4,4,5
ISTA,CUDA,6.0,5.5,6.1,4.7,4.5,5.0,4,4,5
FISTA,OpenMP,2.7,2.6,2.8,3.7,3.6,4.1,5,4,5
FISTA,CUDA,10.9,10.5,11.1,5.7,5.5,6.2,5,5,5
\end{filecontents*}
\csvreader[head to column names, late after line=\\]{./data/codegen_lasso_agent_time.csv}{}{\method & \target & $\tr$ & $\trmin$--$\trmax$ & $\opt$ & $\optmin$--$\optmax$ & $\rounds$ & $\roundsmin$--$\roundsmax$}
        \bottomrule
    \end{tabular}
\end{table}

The agent's reports describe both changes, for the median runs of FISTA on OpenMP and of ISTA on CUDA:

\begin{quote}
\begin{minipage}{\linewidth}
\small\raggedright
\texttt{t0 needs one product with A instead of three. \ldots{} So t0 = act(bias + A$\cdot$y) with y = s1 + g(s1 $-$ s0). This cut traffic per iteration from about 125 MB to 52 MB.}

\smallskip
\texttt{2D-tiled A$\cdot$s with the x slab in shared memory \ldots{} x is now read 20 times (6.4 MB) instead of as 83 MB of sectors. The kernel went from 37 to 41 $\mu$s down to about 27 $\mu$s.}
\end{minipage}
\end{quote}

\subsection{Network optimization}\label{subsec:network}

We consider a minimum cost network flow problem with varying demands~\cite{bertsekas1991linear}.
Each edge has a gain that converts the flow it ships, so the constraint matrix is badly scaled.
The demands $x$ are the parameter.
We solve this linear program with PDHG~\cite{chambolleFirstOrderPrimalDualAlgorithm2011} in the form adapted to linear programs~\cite{applegatePracticalLargeScaleLinear2021}.
We also solve it with a variant that adds momentum in the primal variable, whose extrapolation step \fomtools{} eliminates from the VP.
Practical PDHG solvers equilibrate the constraint matrix before iterating but check termination on the original problem~\cite{applegatePracticalLargeScaleLinear2021}.
We verify both methods with and without Ruiz equilibration~\cite{ruizScalingAlgorithmEquilibrate2001} in the same way.
The equilibrated method iterates on the scaled problem, and the VP measures the residual of the original problem.
The scaling is a constant diagonal change of variables, so it adds no steps or variables to the VP.
Because the residual is unscaled before it enters the objective, the bounds of all variants are directly comparable.
Appendix~\ref{app:experiments_network} gives the formulation, the equilibration, the iterations, and the problem data.

\begin{figure}[ht]
\centering
\begin{tikzpicture}
	\begin{groupplot}[
		group style={group size=2 by 1, horizontal sep=0.6cm, y descriptions at=edge left},
		width=0.5\textwidth,
		height=0.36\textwidth,
		xlabel={$K$},
		xmin=1, xmax=39,
		xtick={10,20,30},
		ymode=log,
		ymin=0.1, ymax=1.1,
		ytick={0.1,0.2,0.5,1},
		log ticks with fixed point,
		grid=major,
		title style={font=\small, yshift=-1ex},
	]
		\nextgroupplot[
			title={PDHG},
			ylabel={primal-dual residual},
			legend to name=networklegend,
			legend columns=2,
			legend cell align=left,
			legend style={font=\small, /tikz/every even column/.append style={column sep=0.4cm}},
		]
		\addplot[color=fomtoolsblue, thick, mark=*, mark size=1.2pt] table [x=K, y=unscaled_mip_bound, col sep=comma] {./data/network_residual.csv};
		\addlegendentry{unscaled (worst-case)}
		\addplot[color=fomtoolsblue, dashed] table [x=K, y=unscaled_prob95, col sep=comma] {./data/network_residual.csv};
		\addlegendentry{unscaled (probabilistic)}
		\addplot[color=red!70!black, thick, mark=*, mark size=1.2pt] table [x=K, y=ruiz_mip_bound, col sep=comma] {./data/network_residual.csv};
		\addlegendentry{Ruiz equilibrated (worst-case)}
		\addplot[color=red!70!black, dashed] table [x=K, y=ruiz_prob95, col sep=comma] {./data/network_residual.csv};
		\addlegendentry{Ruiz equilibrated (probabilistic)}
		\nextgroupplot[title={momentum PDHG}]
		\addplot[color=fomtoolsblue, thick, mark=*, mark size=1.2pt] table [x=K, y=unscaled_momentum_mip_bound, col sep=comma] {./data/network_residual.csv};
		\addplot[color=red!70!black, thick, mark=*, mark size=1.2pt] table [x=K, y=ruiz_momentum_mip_bound, col sep=comma] {./data/network_residual.csv};
		\addplot[color=fomtoolsblue, dashed] table [x=K, y=unscaled_momentum_prob95, col sep=comma] {./data/network_residual.csv};
		\addplot[color=red!70!black, dashed] table [x=K, y=ruiz_momentum_prob95, col sep=comma] {./data/network_residual.csv};
	\end{groupplot}
	\node at ($(group c1r1.north)!0.5!(group c2r1.north) + (0, 1.2cm)$) {\pgfplotslegendfromname{networklegend}};
\end{tikzpicture}
\caption{Worst-case primal-dual residual~\eqref{eq:network_residual} of the original LP~\eqref{prob:network_flow} for PDHG (left) and momentum PDHG (right), without and with Ruiz equilibration, computed by the MIP backend (solid), and probabilistic bounds on the $95$th percentile of the residual over sampled demands (dashed), described in \Sec\ref{subsec:experiments_setup}.}
\label{fig:network_residual}
\end{figure}
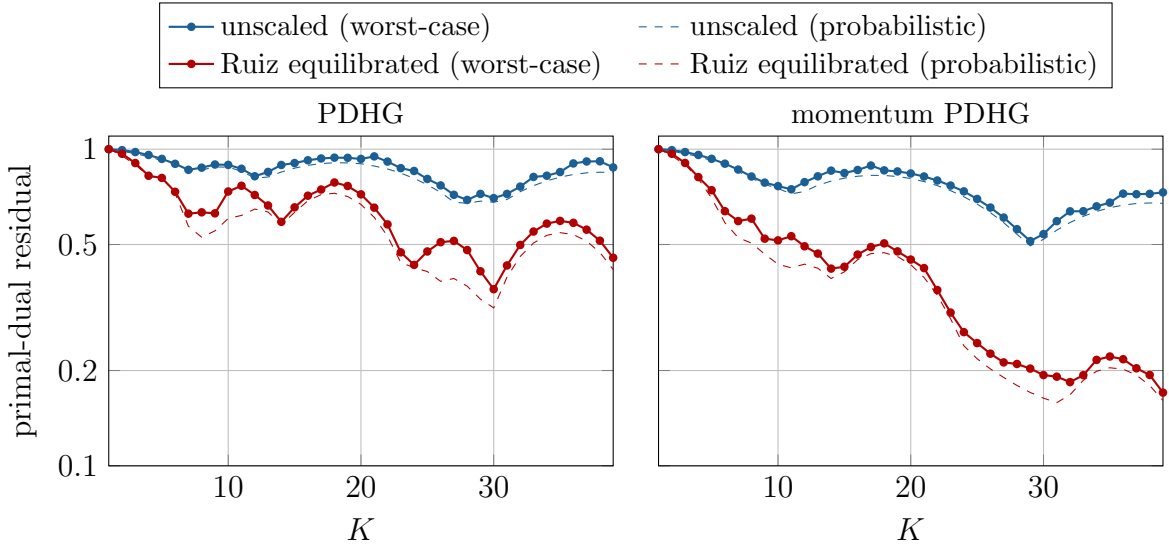

\paragraph{Results.}
Figure~\ref{fig:network_residual} shows that equilibration lowers the worst-case residual of both methods at every $K \ge 2$.
At $K = 39$, it halves the bound of PDHG, from $0.88$ to $0.45$, and lowers that of momentum PDHG by a factor of four, from $0.73$ to $0.17$.
Momentum lowers the worst case at $K = 39$ by $17\%$ without equilibration and by $62\%$ with it.

Table~\ref{tab:codegen_network} shows that, with equilibration, the scheduled momentum coefficient slows the translated solver, as for FISTA.
Unlike for FISTA, the code optimization phase does not close this gap.

\begin{table}[ht]
    \centering
    \caption{Generated C solvers for the verified network instance: time per $200$-iteration solve, and speedups of the translated solver over NumPy, of the optimized solver over the translated one, and in total.
    Each entry is the median over three independent agent runs, and the last column gives the range of the total speedup.}
    \label{tab:codegen_network}
    \codegentable{Network}
\end{table}

\paragraph{Loop versus unrolled solvers.}
The momentum coefficient of momentum PDHG changes at every iteration.
If the user records a run without declaring a schedule, the package therefore generates unrolled code, as described in \Sec\ref{subsec:codegen_ir}.
We compare this unrolled solver with a loop that reads the momentum coefficient from a schedule, for the same equilibrated method.
Both solvers run $K = 50$ iterations on a network with $600$ edges, which Appendix~\ref{app:experiments_network} specifies.
Table~\ref{tab:codegen_network_schedule} shows that the loop solver stores about eighteen times fewer coefficients, and its shared library is about eleven times smaller.
The coefficients of the unrolled solver grow with $K$, while the loop solver only adds one schedule entry per iteration.
The loop solver also accepts an iteration cap at run time, whereas the unrolled solver always runs the recorded $K$ iterations.
After translation, the unrolled solver is slightly faster.
Code optimization then works on one iteration block instead of the full run, and the optimized loop solver is $1.4$ times as fast as the optimized unrolled solver.
Each of its three agent runs is faster than each of the three unrolled ones.
For the loop solver of the median run, the agent reports:

\begin{quote}
\begin{minipage}{\linewidth}
\small\raggedright
\texttt{On rows 620..649, su\_state\_\_d\_xi = -su\_t\_t0\_\_d\_xi, with the same pattern and all columns below 600. So xi*(D x + E t) becomes xi * E (t - x), which removes one product per entry.}
\end{minipage}
\end{quote}

\begin{table}[ht]
    \centering
    \caption{Generated solvers for Ruiz-equilibrated momentum PDHG on a network with $600$ edges, at $K = 50$ iterations per solve in double precision, as a loop that reads the momentum coefficient from a schedule or as unrolled code.
    The steps and coefficient arrays describe one iteration for the loop solver and the full run for the unrolled solver.
    The memory columns give the sizes of the schedule table, the stored coefficients, and the optimized shared library.
    Times are medians over seven interleaved rounds and over three independent agent runs.
    The instance is given in Appendix~\ref{app:experiments_network}.}
    \label{tab:codegen_network_schedule}
    \small
    \setlength{\tabcolsep}{4pt}
    \begin{tabular}{@{}lrrrrrrrrr@{}}
        \toprule
        & & & \multicolumn{3}{c}{Time per solve ($\mu$s)} & & \multicolumn{3}{c}{Memory (KB)}\\
        \cmidrule(lr){4-6}\cmidrule(l){8-10}
        Solver & steps & arrays & NumPy & transl. & optim. & speedup & table & coeff. & library\\
        \midrule
        \begin{filecontents*}{./data/codegen_network_schedule.csv}
method,steps,blocks,numpyus,translatedus,optimizedus,speedtotal,schedkb,coefkb,libkb
loop,2,15,669,124.8,31.0,21.6,0.4,68,203
unrolled,149,217,652,108.7,43.9,14.9,--,1202,2255
\end{filecontents*}
\csvreader[head to column names, late after line=\\]{./data/codegen_network_schedule.csv}{}{\method & \steps & \blocks & \numpyus & \translatedus & \optimizedus & $\speedtotal\times$ & \schedkb & \coefkb & \libkb}
        \bottomrule
    \end{tabular}
\end{table}

\subsection{Optimal control}\label{subsec:opt_control}

We consider the control of a constrained time-invariant linear dynamical system over a finite horizon of $h$ time steps~\cite[Ch.\ 11.3]{borrelliPredictiveControlLinear2017}.
The dynamics are those of a quadcopter linearized around hover~\cite{kouzoupisProperAssessmentQP2015}.
The initial state $x$ is the parameter, with the roll and pitch angles and the altitude varying.
A model predictive controller solves this problem once per sampling period, so the sampling rate fixes the iteration budget $K$.
The question of \Sec\ref{sec:introduction} then becomes a hard real-time requirement.
We solve the problem with the restarted accelerated primal-dual hybrid gradient method (rAPDHG) of~\cite{luPracticalOptimalFirstOrder2026}, with restart periods $K_\mathrm{r} = 6$ and $12$.
Its iterations solve no linear system.
The metric is a residual of the KKT conditions of the problem.
Appendix~\ref{app:experiments_opt_control} gives the formulation, the iterations, the residual map, and the problem data.

\begin{figure}[ht]
\centering
\begin{tikzpicture}
	\begin{semilogyaxis}[
		width=0.64\textwidth,
		height=0.42\textwidth,
		xlabel={$K$},
		ylabel={KKT residual},
		xmin=1, xmax=75,
		grid=major,
		legend cell align=left,
		legend style={font=\small, at={(1.03,1)}, anchor=north west},
	]
		\addplot[color=fomtoolsblue, thick, mark=*, mark size=1.2pt, mark repeat=5, mark phase=5] table [x=K, y=bound_R12, col sep=comma] {./data/control_rapdhg.csv};
		\addlegendentry{$K_\mathrm{r} = 12$ (worst-case)}
		\addplot[color=fomtoolsblue, dashed] table [x=K, y=prob95_R12, col sep=comma] {./data/control_rapdhg.csv};
		\addlegendentry{$K_\mathrm{r} = 12$ (probabilistic)}
		\addplot[color=orange!85!black, thick, mark=*, mark size=1.2pt, mark repeat=5, mark phase=5] table [x=K, y=bound_R6, col sep=comma] {./data/control_rapdhg.csv};
		\addlegendentry{$K_\mathrm{r} = 6$ (worst-case)}
		\addplot[color=orange!85!black, dashed] table [x=K, y=prob95_R6, col sep=comma] {./data/control_rapdhg.csv};
		\addlegendentry{$K_\mathrm{r} = 6$ (probabilistic)}
	\end{semilogyaxis}
\end{tikzpicture}
\caption{Worst-case KKT residual~\eqref{eq:rapdhg_residual} of rAPDHG for the quadcopter controller after $K$ iterations from a cold start, computed by the MIP backend (solid), and probabilistic bounds on the $95$th percentile of the residual over sampled initial states (dashed), described in \Sec\ref{subsec:experiments_setup}.
rAPDHG restarts every $K_\mathrm{r}$ iterations.}
\label{fig:control_rapdhg}
\end{figure}
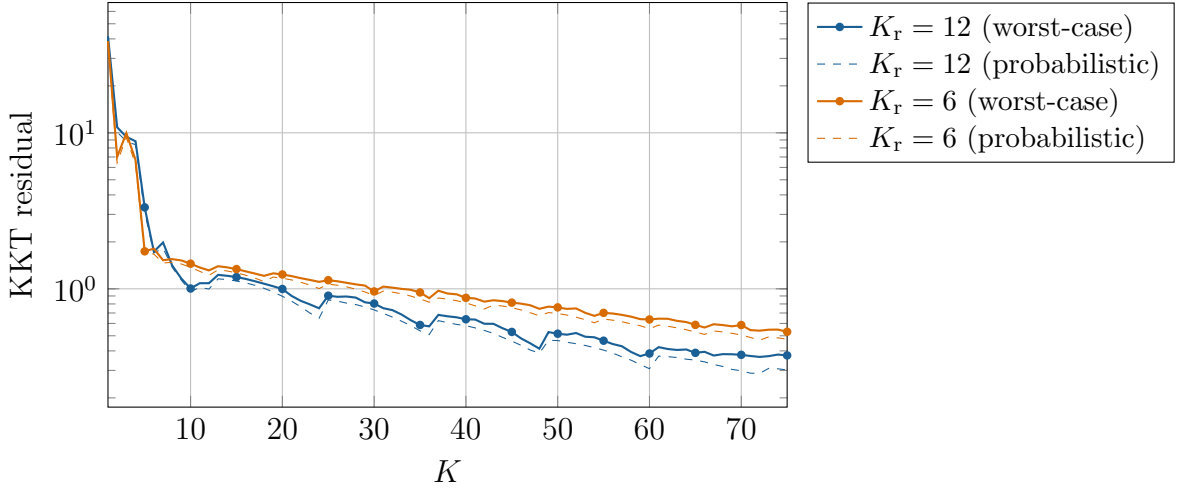

\paragraph{Results.}
Figure~\ref{fig:control_rapdhg} shows that the restart period matters.
From $K = 8$ on, restarting every $K_\mathrm{r} = 12$ iterations gives the smaller worst case at every $K$.
It certifies a residual below $1$ from $K = 20$ on, against $K = 33$ for $K_\mathrm{r} = 6$.
Because the VP bounds every iteration count, it also shows that the worst case is not monotone and rises right after the restarts.

Table~\ref{tab:codegen_control} shows that translation and code optimization each contribute about an order of magnitude.

\begin{table}[ht]
    \centering
    \caption{Generated C solvers for the verified quadcopter instance: time per $200$-iteration solve, and speedups of the translated solver over NumPy, of the optimized solver over the translated one, and in total.
    Each entry is the median over three independent agent runs, and the last column gives the range of the total speedup.}
    \label{tab:codegen_control}
    \codegentable{Control}
\end{table}

\paragraph{Closed-loop flight.}
A real-time controller has a fixed computation budget per sampling instant, so a faster solver buys more iterations.
We generate rAPDHG for a one-second horizon at a $10$\,ms sampling period.
The step sizes and the restart are declared as schedules, so the solver is a single loop.
We fly a simulated quadcopter around a figure eight with $40\,\mu$s of computation per sampling instant, as detailed in Appendix~\ref{app:experiments_opt_control}.
Figure~\ref{fig:control_codegen_flight} shows that the NumPy controller loses the path and the translated solver cuts its corners, while the optimized solver tracks it.
The code optimization phase made each iteration $22$ to $29$ times faster in three independent agent runs, by exploiting the structure of the coefficient blocks, as the report of the run in the figure describes:

\begin{quote}
\begin{minipage}{\linewidth}
\small\raggedright
\texttt{The L block is lower-triangular Toeplitz with 10 distinct coefficients, paired across two input channels. So L$\cdot$w and L$^\mathtt{T}\cdot$v cost 55 two-wide FMAs each, instead of 440 sparse products each. \ldots{}
S is the 40x40 symmetric block of t1\_state\_\_d\_eta. \ldots{} Sum/difference transforms reduce it to 400 FMAs instead of 1600.}
\end{minipage}
\end{quote}

\begin{figure}[ht]
\centering
\includegraphics[width=\textwidth]{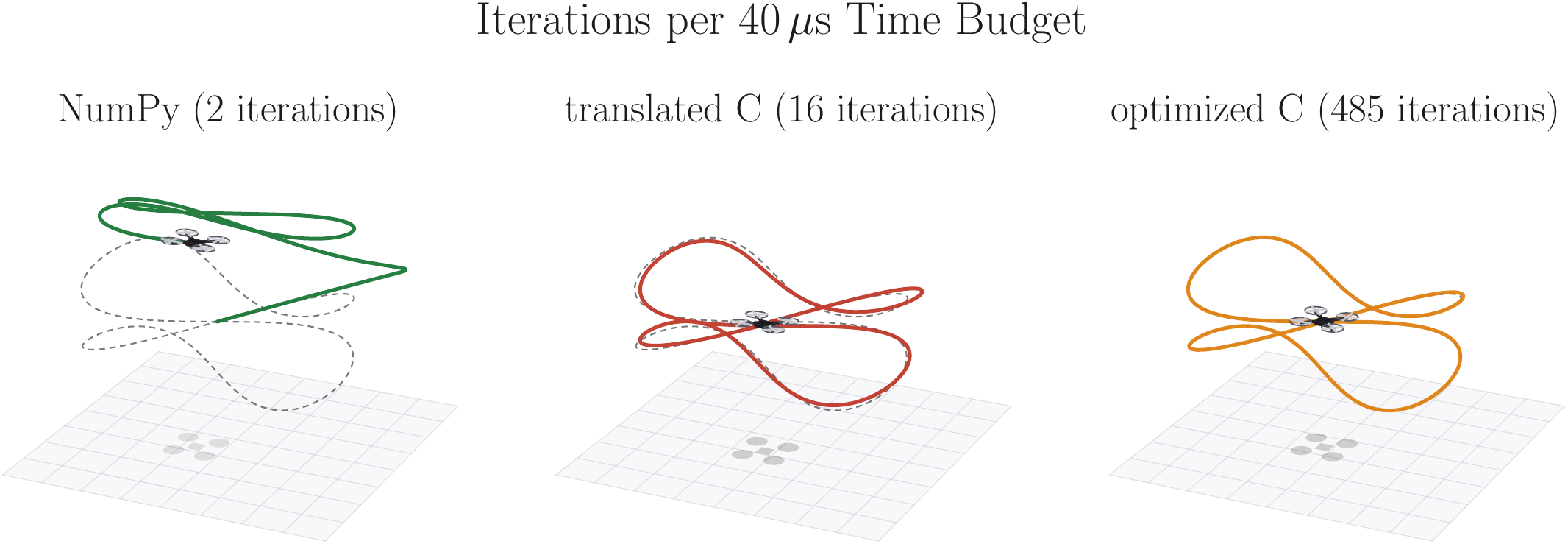}
\caption{Closed-loop flights of the quadcopter on a figure-eight reference (dashed) with the NumPy reference execution, the translated C solver, and the optimized C solver of rAPDHG, each given $40\,\mu$s of computation per sampling instant.
The budget buys $2$, $16$, and $485$ iterations, at roughly $14$, $2.4$, and $0.082\,\mu$s per iteration.
The figure shows the agent run whose optimized solver runs the median number of iterations over three independent runs.
}
\label{fig:control_codegen_flight}
\end{figure}

\subsection{Portfolio rebalancing}\label{subsec:portfolio}

We consider a portfolio rebalancing problem~\cite{markowitz,boydMultiPeriodTradingConvex2017}.
An investor chooses new holdings that trade the risk of the portfolio against its expected return, net of proportional transaction costs, subject to a budget and position limits.
The risk is given by a factor model.
The current holdings and the return forecast are parameters, which change at every solve.
Appendix~\ref{app:experiments_portfolio} gives the formulation and the problem data.

This example tests code generation at scale through the CVXPYgen front end of \Sec\ref{subsec:codegen_cvxpygen}, on problems with up to $6100$ variables.
We compare the \fomtools{} solver of the \texttt{admm} builder with the OSQP solver that CVXPYgen generates.
Both solvers run the same ADMM iteration, so the comparison isolates the implementation.

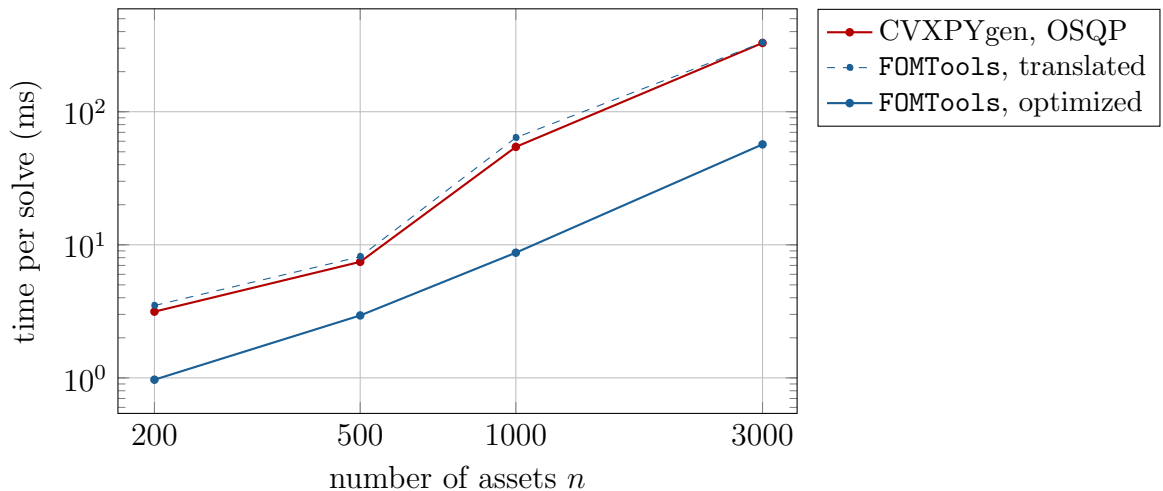
\begin{figure}[ht]
\centering
\begin{tikzpicture}
	\begin{loglogaxis}[
		width=0.64\textwidth,
		height=0.42\textwidth,
		xlabel={number of assets $n$},
		ylabel={time per solve (ms)},
		xmin=170, xmax=3500,
		xtick={200,500,1000,3000},
		xticklabels={$200$,$500$,$1000$,$3000$},
		grid=major,
		legend cell align=left,
		legend style={font=\small, at={(1.03,1)}, anchor=north west},
	]
		\addplot[color=red!70!black, thick, mark=*, mark size=1.2pt] table [x=N, y=osqpms, col sep=comma] {./data/codegen_portfolio_tolerance.csv};
		\addlegendentry{CVXPYgen, OSQP}
		\addplot[color=fomtoolsblue, dashed, mark=*, mark size=1.2pt] table [x=N, y=fomms, col sep=comma] {./data/codegen_portfolio_tolerance.csv};
		\addlegendentry{\fomtools{}, translated}
		\addplot[color=fomtoolsblue, thick, mark=*, mark size=1.2pt] table [x=N, y=fomoptms, col sep=comma] {./data/codegen_portfolio_tolerance.csv};
		\addlegendentry{\fomtools{}, optimized}
	\end{loglogaxis}
\end{tikzpicture}
\caption{Time per solve to tolerance $10^{-4}$ of the generated C solvers for portfolio rebalancing with $n$ assets, on logarithmic axes: the OSQP solver generated by CVXPYgen, and the \fomtools{} solver after translation (dashed) and after code optimization.
All solvers run the same ADMM iteration and stop on comparable residual tests.
Times are medians over three interleaved rounds of the mean time per solve over $100$ draws of the holdings and the return forecast, in double precision, on an Apple M4 Pro, and the generated times are medians over three independent agent runs.}
\label{fig:portfolio_codegen_scaling}
\end{figure}

\paragraph{Results.}
We generate both solvers for $n = 200$, $500$, $1000$, and $3000$ assets.
They use the same constants, start from zero, and stop on comparable residual tests with tolerance $10^{-4}$.
Their mean iteration counts agree within $10\%$, so the comparison measures the cost of an iteration.
Figure~\ref{fig:portfolio_codegen_scaling} shows that the translated solver is $1$ to $17\%$ slower than OSQP, while the optimized solver is $3.3$, $2.5$, $6.2$, and $5.8$ times faster for $n = 200$, $500$, $1000$, and $3000$.
These are medians over three independent agent runs for each size.
The speedups of all twelve runs range from $2.5$ to $7.3$, except in one run for $n = 200$.
There, the rewrite used a division, which the embedded target forbids, so the package kept the translated solver.
Going from the CVXPY problem to the optimized solver takes $8$ to $43$ minutes.
On the larger instances, most of the solve time goes to the triangular solves with the factorization of the linear system in~\eqref{eq:osqp_iterations}.
The agent adapted the indexing of these solves to the structure of the factor, as its report for the median run with $n = 1000$ describes:

\begin{quote}
\begin{minipage}{\linewidth}
\small\raggedright
\texttt{Dense column tails in the LDL\^{}T solve (59.0 $\to$ 22.5 ms). 308k of the factor's 315k nonzeros sit in contiguous row runs at the end of each column, mostly rows 6100..6200. \ldots{} The forward pass does each run as a contiguous axpy and the backward pass does it as a dot product with 8 accumulators, so both vectorize.}
\end{minipage}
\end{quote}

\section{Conclusion}\label{sec:conclusion}

We presented \fomtools{}, a package for developing, analyzing, and deploying first-order methods from one mathematical specification.
Its shared affine and PWA representation supports exact MIP formulations for worst-case analysis, batched simulation with probabilistic bounds, and specialized code generation.
A coding agent translates the iteration and then optimizes its implementation, and the package runs the code after both phases and compares it with a reference on sampled inputs.
In our experiments, the verification bounds revealed worst-case behavior that sampling misses.
The generated solvers were up to two orders of magnitude faster than NumPy implementations and outperformed hand-written GPU code and the OSQP solvers generated by CVXPYgen.

MIP verification remains costly as the horizon and problem size grow.
Specialized branch-and-bound procedures could exploit the structure of the algorithm iterations to reduce this cost.
Code generation requires fixed matrices, and its numerical checks are finite tests.
Connecting the VP bounds to floating-point executions of the generated code would require a separate error analysis.
Supporting further operations would extend both analysis and code generation to a broader class of first-order methods.

\section*{Acknowledgements}
This work was supported in part by the National Science Foundation under CAREER award ECCS-2239771, the Office of Naval Research under Young Investigator Program award N00014-25-1-2147, and the Alfred P. Sloan Foundation through a Sloan Research Fellowship under grant FG-2026-79479.

\appendix

\section{Implementation notes}\label{app:implementation_notes}

The expression and step representations are shared by the analysis backends and code generation.
Each analysis backend implements an interface for variables and bounds, constraints, objectives, and results.
A solver wrapper translates the resulting optimization model into the solver's native API.
For simulation, each step supplies a numerical evaluation rule.
Coefficient matrices can be dense arrays or SciPy sparse matrices.
Expression bounds are recomputed when a backend uses them, so tighter bounds propagate through later affine operations.
After a MIP solve, expressions can also be evaluated at a feasible incumbent, including their constant terms and affine combinations.

\section{Experimental SDP backend}\label{app:sdp_backend}

\fomtools{} also ships an experimental backend that forms the semidefinite relaxation of~\cite{ranjanVerificationFirstorderMethods2025}.
The user selects it with \texttt{method="sdp"}, writes the same code as for the other backends, and solves the resulting SDP relaxation with a solver such as Clarabel.
Let $y$ stack every leaf variable of the VP.
The relaxation replaces the rank-one matrix $(1, y)(1, y)^\tpose$ with a positive semidefinite matrix $X$ that satisfies $X_{00} = 1$.
The step constraints and the bounds become linear in $X$, so the optimal value of the relaxation is an upper bound on the optimal value of~\eqref{eq:vp}.

The backend natively supports only ReLU steps, and it rewrites the other PWA steps into ReLU primitives through
\begin{equation}\label{eq:step_decompositions}
	\begin{array}{rl}
		\softthresh_\lambda(v) &= \relu(v - \lambda) - \relu(-v - \lambda),\\
		\sat_{[l,u]}(v) &= \relu(v - l) - \relu(v - u) + l.
	\end{array}
\end{equation}
The SDP directly handles the squared $\ell_2$-norm objective.
It handles the $\ell_1$ and $\ell_\infty$ objectives the same way as~\eqref{eq:step_decompositions}, with $2n$ and $3n - 1$ auxiliary ReLU steps for an $n$-dimensional expression.
It includes the bound products, the triangle inequalities, and the linear ReLU envelope of~\cite{ranjanVerificationFirstorderMethods2025}, together with products of general inequality constraints and optional interpolation inequalities between steps.
The OBBT bounds of \Sec\ref{subsec:mip_backend} tighten all of them.
Each constraint involves the principal submatrix of $X$ on the indices of its expressions.
By default, the canonicalizer imposes positive semidefiniteness only on these principal submatrices rather than on all of $X$.
This gives a smaller but possibly weaker relaxation.
The option \texttt{full\_psd\_cone=True} imposes it on the full matrix.

The SDP backend remains experimental because numerical solves can be unreliable even for moderately sized problems.
Its semidefinite cones often share most of their indices, and the resulting formulation can be degenerate.
Improving the formulation's conditioning and integrating more reliable SDP solvers remain future work.

\section{Proximal API reference}\label{app:prox_api}

This appendix lists every proxable function that \fomtools{} ships with.
Throughout, $v$ is the input expression, $u$ is the output, and $\gamma > 0$ is the proximal step size of~\eqref{eq:prox_def}.
We write $I_C$ for the indicator of a set $C$, which is $0$ on $C$ and $+\infty$ elsewhere, and $\Pi_C$ for the Euclidean projection onto $C$.
The symbol $\odot$ denotes the element-wise product, and $\perp$ denotes complementarity, so that $0 \le a \perp c \ge 0$ means $a \ge 0$, $c \ge 0$, and $a \odot c = 0$.

Table~\ref{tab:prox_closed_form} lists the functions whose proximal operator has a closed form in the three PWA primitives of \Sec\ref{subsec:steps}.
Table~\ref{tab:prox_kkt} lists the remaining functions, which the package represents through a linear system or the optimality conditions of the proximal problem.
For these the package introduces the primal and dual variables of that problem and encodes its optimality constraints within the VP constraints.
Table~\ref{tab:prox_sets} gives the sets accepted by \texttt{project}, each of which is a thin wrapper over the corresponding indicator.

\begin{table}[ht]
    \centering
    \caption{Proxable functions with a closed-form proximal operator in the ReLU, soft-thresholding, and saturated linear primitives.
    The saturated linear step $\sat_{[l,u]}$ clips element-wise to $[l, u]$, and $\softthresh_\lambda$ is element-wise soft thresholding.
    Bracketed constructor arguments are the Python parameter names.
    Weights, radii, and regularization parameters are nonnegative, $l \le u$, $y_i \in \{-1, 1\}$, $0 < \alpha \le 1$, and $a \ne 0$.}
    \label{tab:prox_closed_form}
    \small
    \begin{tabular}{@{}lll@{}}
        \toprule
        Proxable function & $f(x)$ & $\prox_{\gamma f}(v)$\\
        \midrule
        \multicolumn{3}{@{}l@{}}{\itshape No step: affine arithmetic only}\\
        \texttt{Linear(c)} & $c^\tpose x$ & $v - \gamma c$\\
        \texttt{SqL2Norm()} & $(1/2)\|x\|_2^2$ & $v / (1 + \gamma)$\\
        \addlinespace
        \multicolumn{3}{@{}l@{}}{\itshape One step}\\
        \texttt{IndNonNegative()} & $I_{\{x \ge 0\}}$ & $\relu(v)$\\
        \texttt{IndNonPositive()} & $I_{\{x \le 0\}}$ & $-\relu(-v)$\\
        \texttt{IndBox(l, u)} & $I_{\{l \le x \le u\}}$ & $\sat_{[l,u]}(v)$\\
        \texttt{IndBallLinf(r)} & $I_{\{\|x\|_\infty \le r\}}$ & $\sat_{[-r,r]}(v)$\\
        \texttt{IndPoint(p)} & $I_{\{x = p\}}$ & $\sat_{[p,p]}(v) = p$\\
        \texttt{NormL1(lambd)} & $\sum_i \lambda_i |x_i|$ & $\softthresh_{\gamma\lambda}(v)$\\
        \texttt{IndHalfSpace(a, b)} & $I_{\{a^\tpose x \le b\}}$
            & $v - \relu\!\left((a^\tpose v - b) / \|a\|_2^2\right) a$\\
        \addlinespace
        \multicolumn{3}{@{}l@{}}{\itshape One step plus affine arithmetic}\\
        \texttt{ElasticNet(alpha, beta)} & $\alpha\|x\|_1 + \beta\|x\|_2^2$
            & $\softthresh_{\alpha\gamma}(v) / (1 + 2\beta\gamma)$\\
        \texttt{HuberLoss(delta)} & $\sum_i h_\delta(x_i)$
            & $\left(v + \gamma\,\softthresh_{\delta(1+\gamma)}(v)\right) / (1 + \gamma)$\\
        \texttt{SoftNonNegative(alpha)} & $((1-\alpha)/2\alpha)\sum_i \relu(-x_i)^2$
            & $\tilde\alpha\, v + (1 - \tilde\alpha)\,\relu(v)$\\
        \texttt{HingeLoss()} & $\sum_i \relu(1 - x_i)$
            & $v + \sat_{[0,\gamma]}(1 - v)$\\
        \texttt{LabeledHingeLoss(y)} & $\sum_i \relu(1 - y_i x_i)$
            & $v + y \odot \sat_{[0,\gamma]}(1 - y \odot v)$\\
        \texttt{SumPositive()} & $\sum_i \relu(x_i)$
            & $v - \sat_{[0,\gamma]}(v)$\\
        \texttt{SumNegative()} & $\sum_i \relu(-x_i)$
            & $v - \sat_{[-\gamma,0]}(v)$\\
        \bottomrule
    \end{tabular}
\end{table}

The two entries of Table~\ref{tab:prox_closed_form} that need a definition are the Huber function~\cite{huber1964huber}
\begin{equation*}
    h_\delta(t) = \begin{cases}
        (1/2) t^2 & |t| \le \delta\\
        \delta |t| - (1/2)\delta^2 & |t| > \delta,
    \end{cases}
\end{equation*}
and the effective slope $\tilde\alpha = \alpha / (\alpha + \gamma(1 - \alpha))$ of \texttt{SoftNonNegative}, whose proximal operator is a leaky ReLU.
Two limits of that slope are worth noting: at $\gamma = 1$ it is the standard leaky ReLU with slope $\alpha$, and as $\alpha \to 0^+$ the penalty on the negative orthant grows without bound and the operator recovers the plain ReLU of \texttt{IndNonNegative}.

\begin{table}[ht]
    \centering
    \caption{Proxable functions encoded through the optimality conditions of~\eqref{eq:prox_def} rather than a PWA step.
    The output is the primal variable $u$; the remaining symbols are auxiliary variables that the package introduces.
    The matrix $Q$ is positive semidefinite.
    A dagger marks the one function for which the user must supply a bound on the multipliers.
    The bound must keep at least one valid multiplier for every proximal input.}
    \label{tab:prox_kkt}
    \small
    \begin{tabular}{@{}l l >{\raggedright\arraybackslash}p{0.44\textwidth}@{}}
        \toprule
        Proxable function & $f(x)$ & Encoding\\
        \midrule
        \texttt{Quadratic(Q, q)} & $(1/2) x^\tpose Q x + q^\tpose x$
            & $(I + \gamma Q)\, u = v - \gamma q$\\
        \addlinespace
        \texttt{LeastSquares(A, b)} & $(1/2)\|Ax - b\|_2^2$
            & As \texttt{Quadratic} with $Q = A^\tpose A$ and $q = -A^\tpose b$\\
        \addlinespace
        \texttt{IndAffine(A, b)} & $I_{\{Ax = b\}}$
            & $u + A^\tpose \nu = v$, \quad $A u = b$\\
        \addlinespace
        \texttt{IndPolyhedral(A, b)}$^\dagger$ & $I_{\{Ax \le b\}}$
            & $u = v - A^\tpose \lambda$, \quad $0 \le \lambda \perp b - A u \ge 0$\\
        \addlinespace
        \texttt{IndSimplex(a)} & $I_{\{x \ge 0,\; \ones^\tpose x = a\}}$
            & $u = v + \lambda - \nu \ones$, \quad $\ones^\tpose u = a$, \quad $0 \le \lambda \perp u \ge 0$; the multipliers are bounded by the input box\\
        \addlinespace
        \texttt{IndBallL1(r)} & $I_{\{\|x\|_1 \le r\}}$
            & $u = v - \lambda_1 + \lambda_2$, \quad $\lambda_1 + \lambda_2 = \lambda_3 \ones$, \quad
              $0 \le (\lambda_1, \lambda_2, \lambda_3) \perp (t - u,\; t + u,\; r - \ones^\tpose t) \ge 0$\\
        \addlinespace
        \texttt{NormInf(lambd)} & $\lambda \|x\|_\infty$
            & $u = v - \Pi_{\{\|\cdot\|_1 \le \gamma\lambda\}}(v)$ by Moreau decomposition,
              reusing \texttt{IndBallL1}\\
        \addlinespace
        \texttt{TotalVariation1D(lambd)} & $\lambda \|Dx\|_1$
            & $u = v - D^\tpose w$, \quad $D u = z^+ - z^-$, \quad
              $0 \le \gamma\lambda \mp w \perp z^\pm \ge 0$\\
        \bottomrule
    \end{tabular}
\end{table}

The auxiliary variables of Table~\ref{tab:prox_kkt} are the multipliers of the corresponding constraints: $\nu$ for an equality, $\lambda$ for an inequality, and $\lambda_1$, $\lambda_2$, $\lambda_3$ for the three groups of inequalities of the $\ell_1$-ball.
The $\ell_1$-ball formulation lifts the problem to the pair $(u, t)$, where $t$ bounds $|u|$ element-wise, which is what makes the constraint set polyhedral.
For \texttt{TotalVariation1D}, the matrix $D$ is the first-difference operator, $w$ is the multiplier of $Du = z^+ - z^-$, and the split $z^+, z^- \ge 0$ carries the two branches of the absolute value.
The nonnegativity halves of its complementarity conditions encode the dual box $-\gamma\lambda \le w \le \gamma\lambda$, so all of its bounds follow from the problem data.
Bounds on the multipliers matter because the MIP backend encodes each complementarity condition with a big-$M$ constraint, so a loose dual bound weakens the formulation exactly as a loose iterate bound does in \Sec\ref{subsec:mip_backend}.
The package derives these bounds from the input box wherever the geometry allows it, and OBBT tightens them further.

\begin{table}[ht]
    \centering
    \caption{Sets accepted by \texttt{project}.
    Each is a wrapper that returns the indicator in the last column, so \texttt{project(C, expr)} and \texttt{prox(C.indicator(), expr)} are the same call.
    The parameterless sets may also be passed as a bare class.}
    \label{tab:prox_sets}
    \small
    \begin{tabular}{@{}lll@{}}
        \toprule
        Set & $C$ & Indicator\\
        \midrule
        \texttt{Box(l, u)} & $\{ x \mid l \le x \le u \}$ & \texttt{IndBox}\\
        \texttt{NonNegative()} & $\{ x \mid x \ge 0 \}$ & \texttt{IndNonNegative}\\
        \texttt{NonPositive()} & $\{ x \mid x \le 0 \}$ & \texttt{IndNonPositive}\\
        \texttt{Point(p)} & $\{ p \}$ & \texttt{IndPoint}\\
        \texttt{LinfBall(r)} & $\{ x \mid \|x\|_\infty \le r \}$ & \texttt{IndBallLinf}\\
        \texttt{L1Ball(r)} & $\{ x \mid \|x\|_1 \le r \}$ & \texttt{IndBallL1}\\
        \texttt{HalfSpace(a, b)} & $\{ x \mid a^\tpose x \le b \}$ & \texttt{IndHalfSpace}\\
        \texttt{Polyhedron(A, b)} & $\{ x \mid Ax \le b \}$ & \texttt{IndPolyhedral}\\
        \texttt{Simplex(a)} & $\{ x \mid x \ge 0,\; \ones^\tpose x = a \}$ & \texttt{IndSimplex}\\
        \texttt{AffineSubspace(A, b)} & $\{ x \mid Ax = b \}$ & \texttt{IndAffine}\\
        \bottomrule
    \end{tabular}
\end{table}

Every function in both tables is available on every backend.
The MIP backend encodes the tables as written, while the experimental SDP backend first rewrites soft-thresholding and saturated linear steps into ReLU primitives (Appendix~\ref{app:sdp_backend}).
The simulator instead evaluates each function numerically per sample, so the entries of Table~\ref{tab:prox_closed_form} become element-wise array operations and those of Table~\ref{tab:prox_kkt} become a batched subproblem solve.
One restriction applies to the simulator: the scalar arguments \texttt{r}, \texttt{a}, and \texttt{lambd} of the functions in Table~\ref{tab:prox_kkt} must be constants rather than problem parameters, because the batched evaluation fixes them across the sample.

\section{Code generation details}\label{app:codegen_files}

\subsection{Algorithm specification}\label{app:codegen_inputs}

Code generation accepts three inputs: a recorded run, a single stationary iteration, or one iteration with scheduled coefficients.
\fomtools{} uses the variable declarations to distinguish the state carried between iterations from the parameters supplied at solve time.
The user can also specify the state and parameters explicitly through \texttt{state} and \texttt{params}.

\paragraph{A recorded run.}
For the FISTA example in \Sec\ref{subsec:workflow}, a single call generates a solver from the recorded iterates.

\codelisting[language=fomPython]{code/codegen.py}

A recorded run stores every iteration separately, even when all iterations perform the same update.
Before generating code, the package checks whether they do.
It applies the first recorded iteration repeatedly and compares the result with the recorded iterates on sampled parameter values, within a numerical tolerance.
If they agree, the solver implements the method as a loop over a single copy of the iteration.
Otherwise, the package generates unrolled code, which writes out each recorded iteration in sequence with its own coefficients.
In both cases, the solver executes the recorded number of iterations by default.
The user can request unrolled code directly with \texttt{roll=False}.
Unrolled code stores coefficient arrays that are identical across iterations only once.

\paragraph{A stationary iteration.}
For a stationary method, the user records a single iteration and calls \texttt{codegen\_rolled}, as for ISTA in \Sec\ref{subsec:codegen_ir}.
If the user also passes a recorded trajectory through \texttt{iterates} and the parameter values through \texttt{param\_values}, the package checks, at those parameter values, that repeating the iteration reproduces the trajectory.

\paragraph{Scheduled coefficients.}
For FISTA, the user writes one iteration as a function of the momentum coefficient $\gamma$ and returns a \texttt{Recording} with the current and next states.
The schedule gives $\gamma_k$ as a function of the zero-based iteration index $k$ or as a sequence of values.
The following example uses the momentum sequence of~\eqref{eq:fista_iterations} and carries the two latest iterates in its state.

\codelisting[language=fomPython]{code/codegen_schedule.py}

To find $M_0$ and $M_1$ in each coefficient matrix $M(\gamma)=M_0+\gamma M_1$, the package records the iteration at two values of $\gamma$ and checks the affine fit at a third.
Each additional scheduled scalar requires one more recording.
Only the coefficient matrices may depend on the schedule.
The step types, dimensions, PWA thresholds and bounds, matrices of linear solves, and constant offsets must not.
At generation time, the package stores the schedule values in a table up to the iteration cap.
If the cap is raised at run time beyond the table, the solver holds the last tabulated value.

\subsection{Routines, targets, and interface}

\paragraph{Supporting routines and data.}
For a step that solves $W t=r$ with a fixed matrix $W$, the package factorizes $W$ once and shares the factorization across steps with the same matrix.
It never forms the inverse of $W$.
By default, symmetric systems use $\ldlt$ when the factorization passes numerical checks; otherwise, the package uses LU with partial pivoting.
On the CPU targets, the factors are stored in a sparse format with reciprocal diagonal entries, so the generated substitutions require no division.
During code optimization, the agent may change how these factors are stored and applied, while the reference solutions remain fixed.
The generated matrix-vector products read fixed coefficient arrays.
The thresholds and bounds of PWA steps are also fixed, but bounds that depend on the parameter with a constant width $u-l$ can be handled through the identity $\sat_{[l,u]}(y)=l+\sat_{[0,u-l]}(y-l)$, which moves the parameter into the argument of the step.

\paragraph{Targets.}
The default target is embedded C, which follows the OSQP code generator~\cite{banjacEmbeddedCodeGeneration2017} in using a statically sized workspace, no dynamic allocation, and no division in the iteration.
The agent receives these requirements as instructions.
During code optimization, a static check also rejects allocation and division in rewritten linear-solve and matrix-vector files.
The solver file is held to these requirements by the instructions only.
The OpenMP target runs on a multicore CPU and marks the data-parallel loops with pragmas.
Its build files compile them serially by default, so the checks run on any machine, and \texttt{OMPFLAGS=-fopenmp} enables threads.
The CUDA target runs the whole iteration on an NVIDIA GPU.
Matrix-vector products, linear solves, and reductions all run on the device.
The host reads back at most one scalar per iteration for the stopping test, plus one scalar per CG iteration when a CG solve is present.
Direct solves use the cuDSS sparse direct solver~\cite{nvidiaCuDSS2026}, which factorizes $W$ on the device at the first solve.
For symmetric positive definite systems, both parallel targets can replace the factorization with a Jacobi-preconditioned conjugate gradient (CG) method.
CG stops when the relative residual falls below $10^{-12}$ in double precision and $10^{-6}$ in single precision, or after $2n + 50$ iterations.
The reference execution uses the same CG, so the numerical check compares against this approximate solve.
With GPU execution available, code optimization may rewrite the CUDA solver, matrix-vector products, and linear solves in their \texttt{.cu} files.
The matrix data in \texttt{include/fom\_ops\_data.cuh} and the metric evaluation in \texttt{src/fom\_objective.cu} remain fixed.
Without GPU execution, only the solver file is editable and the check only compiles the code.

\paragraph{Interface.}
Every target exposes the five functions below, which act on one workspace struct.
By default, the solver runs to an iteration cap.
When generating a loop, the user can instead stop early once the infinity norm of the state change, or the verification metric, falls below a tolerance.
For the verification metric, the generated code evaluates the same expression as the verifier.
For a squared $\ell_2$ metric, the tolerance therefore bounds the squared norm.
The tolerance and the interval between metric evaluations can be changed at run time.
Unrolled code always executes the full recorded run.
A parameter update restores the specified initial state, so successive solves do not warm start.
The package generates double-precision code by default and single-precision code on request.
The CVXPY interface supports only double precision.
The generated Makefiles compile with \texttt{-O3}, except for the embedded C library, which uses \texttt{-O2}.
The timing benchmark always uses \texttt{-O3}, and the user can override both flags.

\codelisting[language=C]{code/fom_abi.c}

\paragraph{Checks of rewritten routines.}
When the matrix-vector products are editable, the package tests each routine against its reference product.
It applies the routine to every unit vector when the number of columns times the sum of the nonzeros and rows is at most $10^8$, and to two standard normal vectors in every case.
The output buffer starts nonzero, which checks that the routine adds the product to the existing output.
On the CPU targets, \texttt{fom\_ops.c} also holds one routine per output block that sums its per-iteration products, one that also applies the bias and the PWA step, and one routine for the whole iteration.
The package tests each against the corresponding combination of reference products, on unit vectors at every schedule setting and on random inputs.
Each editable direct linear solve is tested on three standard normal right-hand sides.
CG is checked only through the complete solver.
For a matrix-vector routine, the maximum absolute error must be at most $\epsilon_\mathrm{kern}(1+\|v_\mathrm{ref}\|_\infty)$, where $v_\mathrm{ref}$ is the reference product.
A rewrite of a product may only reorder arithmetic, so $\epsilon_\mathrm{kern}$ is much smaller than $\epsilon_\mathrm{gen}$.
For a direct linear solve, the error must be at most $\epsilon_\mathrm{gen}(1+\|v_\mathrm{ref}\|_\infty)$, where $v_\mathrm{ref}$ is the reference solution.

Source checks keep every coefficient inside these tested routines.
The solver must call the matrix-vector routines instead of reading their coefficient arrays, and each routine and its helpers may read only the matrix data of the products it covers.
The one exception is a product by a $0/1$ selection matrix, which the solver may apply as a copy through a fixed index table.
Such copies are checked only through the complete solver.
New floating-point literals in the solver and in the matrix-vector routines are limited to $0$, $1/2$, $1$, and $2$.
Routines probed on every unit vector are exempt, because those probes recover the matrix they apply entry by entry.
Setup and parameter updates may not add calls to any of these routines, so the agent cannot rebuild a matrix by probing it.
These checks are syntactic.
They complement the numerical tests and do not establish equivalence for all inputs.
After each successful phase, the package writes production copies of the editable files without profiling markers.
Code optimization also writes a report of the rounds, changes, timings, and checks, and of the optimizations the agent considered but the rules forbade.

\paragraph{Default constants.}
The following defaults apply to the code generation procedure in \Sec\ref{sec:codegen}.
\begin{itemize}
    \item \emph{Loop check.} For a recorded run, the test uses three parameter draws, uniform on the parameter bounds with $\pm 1$ in place of an infinite bound, and an absolute tolerance of $10^{-9}$.
    \item \emph{Schedules.} The affine fit must hold to a relative tolerance of $10^{-10}$.
    \item \emph{Factorizations.} The package treats $W$ as symmetric when $\max_{ij} |W_{ij} - W_{ji}| \le 10^{-9} \max\{1, \max_{ij} |W_{ij}|\}$, and then factorizes its symmetric part.
    QDLDL, the factorization used by OSQP~\cite{stellatoOSQPOperatorSplitting2020}, computes $\ldlt$ in sparse arithmetic without pivoting, after an approximate minimum degree permutation $P$.
    The package accepts it only if the backward error $\|L D L^\tpose - P W P^\tpose\|_F / \max\{1, \|W\|_F\}$ is at most $10^{-9}$ and no pivot is smaller than $10^{-12}$ in magnitude.
    Otherwise it falls back to LU with partial pivoting, which is computed in dense arithmetic and then stored sparse, and it rejects $W$ as singular if a diagonal entry of $U$ is at most $10^{-12} \max\{1, \max_i |U_{ii}|\}$ in magnitude.
    On the CUDA target, cuDSS computes the direct factorization on the device instead.
    \item \emph{Matrix-vector products.} A product is generated in compressed sparse row format when the matrix has at least $16$ entries and fewer than half of them are nonzero, and in dense format otherwise.
    \item \emph{CUDA.} The target selects CG for symmetric positive definite systems of dimension at least $256$ and cuDSS, version $0.8$ or later, for the other systems.
    The package links cuDSS only when a direct solve is present.
    The agent writes kernels for the PWA steps, and the package supplies the verification metric.
    \item \emph{Agent.} Translation allows five attempts within a default time budget of $10$ minutes and fails if no attempt passes the check.
    Code optimization runs at most five rounds.
    Each round is a new agent session with up to five attempts, a default budget of $2$ hours, and $200$ turns per attempt for Claude, because builds and benchmarks on a GPU are slow.
    All of these limits can be set for each run.
    The phase stops after a round that gains less than $2\%$.
    Within a round, the agent is instructed to stop after two changes with less than about $2\%$ improvement.
    With profiling enabled, gate~2 rejects a first round that is more than $2\%$ slower than the translated solver, allowing small timing variations, and a later round that is not faster than the best earlier one.
    Claude and Codex run through their software development kits, with their command-line interfaces as a fallback, and pi runs through its remote procedure call interface.
    Any other command-line agent is added through a declarative entry in a configuration file, and other interfaces through a driver class.
    A stub driver supplies a fixed solver loop, so the pipeline also runs offline without a model.
    \item \emph{Numerical acceptance check.} The check uses eight parameter draws and $\epsilon_\mathrm{gen} = 10^{-6}$ in double precision and $10^{-4}$ in single precision.
    The routine tests use $\epsilon_\mathrm{kern} = 10^{-10}$ in double precision and $10^{-4}$ in single precision.
    Infinite parameter bounds are replaced by $\pm 1$ for sampling.
    For unrolled code, it also runs the solver with a larger iteration cap, to confirm that the solver ignores the cap.
    With metric-based stopping, the iteration count and status must match the reference exactly, and the metric must agree within $\epsilon_\mathrm{gen}\max\{1,|\phi_\mathrm{ref}|\}$.
    Additional runs vary the check interval and iteration cap.
    \item \emph{Stopping.} The default iteration cap is $1000$ unless a recorded run supplies its horizon.
    When stopping on the verification metric, the default check interval is $25$ iterations, and the metric is also evaluated at the last iteration.
    \item \emph{CVXPY builders.} The \texttt{admm} builder uses OSQP's defaults: $\rho = 0.1$, $\sigma = 10^{-6}$, over-relaxation $1.6$, and $\rho$ scaled by $10^3$ on equality rows.
    The \texttt{admm} and \texttt{pdhg} builders record three iterations to check the loop.
    The default diagonal scaling uses ten Ruiz equilibration passes and remains fixed across solves.
    Residual-based stopping uses the original problem units, whereas state-change stopping uses the scaled state.
    A parameter without supplied bounds receives its current value plus or minus $1$, with a warning.
\end{itemize}

\section{Problem formulations and experimental details}\label{app:experiments}

This appendix gives the timing protocol of the generated solvers.
For each family of \Sec\ref{sec:experiments}, it also gives the parametric problem, the iterations analyzed, and the problem data.

\paragraph{Timing protocol.}
This protocol applies to the solvers generated for the verified instances, in Tables~\ref{tab:codegen_lasso}, \ref{tab:codegen_network}, and~\ref{tab:codegen_control}.
The larger code generation experiments use different instances, which the last paragraph of each subsection below gives, together with any change to this protocol.
The NumPy baselines are direct single-instance implementations of the same iteration without residual computations, timed as the median over $50$ solves after $3$ warm-up solves at one verified parameter.
The generated code is timed as the mean over up to $5000$ solves that cycle through the parameters of the numerical acceptance check, stopping after $2$ seconds of solve time.
Each reported time is the median over seven interleaved rounds of the three measurements, and across rounds the NumPy, translated, and optimized times varied by at most $14\%$, $12\%$, and $34\%$ of their medians.
Each solver is generated three times, in independent agent runs, and the tables report the medians over the three runs.
Every generated solver in \Sec\ref{sec:experiments} passed the numerical acceptance check of \Sec\ref{subsec:codegen_verify}.
Over $20$ draws from the verified parameter set at $200$ iterations, the largest absolute deviation from the NumPy implementation was $3.0 \cdot 10^{-10}$ for FISTA and at most $1.1 \cdot 10^{-13}$ for every other solver.
All generation and timing ran on an Apple M4 Pro Mac running macOS 26, with Apple clang 21.0.0 and NumPy 2.5.3.
The coding agent was Claude Opus 5.5, run through the Claude software development kit with a budget of $10$ minutes for translation and for each round of code optimization.
Measured over each phase, including the package's checks and the timing comparisons at gate~2, translation took $23$ to $64$ seconds and code optimization $2.2$ to $9.6$ minutes, in two to five rounds.

\subsection{Lasso}\label{app:experiments_lasso}

The Lasso problem~\eqref{prob:lasso} and the ISTA and FISTA iterations are given in \Sec\ref{subsec:lasso_example}.

\paragraph{Algorithms.}
We analyze ISTA~\eqref{eq:ista_iterations} and FISTA~\eqref{eq:fista_iterations}, both with step size $\eta = 1/L$, where $L = \lambdamax(D^\tpose D)$.
Each ISTA iteration is a single soft-thresholding step.
FISTA adds one affine momentum step, which \fomtools{} eliminates to produce the reduced encoding~\eqref{eq:fista_compact_rewrite} with no manual algebra from the user.

\paragraph{Residual.}
Write $w^{K-1}$ for the point from which ISTA or FISTA computes $s^K$, so $w^{K-1} = s^{K-1}$ for ISTA.
The optimality condition of the soft-thresholding step gives $(w^{K-1} - s^K)/\eta - D^\tpose(D w^{K-1} - x) \in \lambda \partial \|s^K\|_1$, and adding $D^\tpose(D s^K - x)$ to both sides shows that
\begin{equation}\label{eq:lasso_residual}
    r^K = (\eta^{-1} I - D^\tpose D)(w^{K-1} - s^K) \in D^\tpose(D s^K - x) + \lambda \partial \|s^K\|_1.
\end{equation}
If $r^K = 0$, then $s^K$ solves~\eqref{prob:lasso}.
The residual is a linear expression of the last two iterates of ISTA and of the last three iterates of FISTA.
We report its worst case $\|r^K\|_\infty$ for both methods.

\paragraph{Verified instance.}
We choose $m = 25$, $n = 50$, and $\lambda = 0.001$, so the dictionary is underdetermined and $D^\tpose D$ is singular.
The dictionary $D$ is sparse.
Each entry is drawn i.i.d.\ from the Gaussian distribution $\mathcal{N}(0, 1/m)$ with probability $0.2$ and is zero otherwise.
We then add one nonzero entry, in a random position, to every column and every row that has none, and normalize the columns, which gives the step size $\eta = 1/L = 0.221$.
Similarly to~\cite[Sec.\ 4.1]{chenTheoreticalLinearConvergence2018}, we generate $100$ sparse vectors $z^\mathrm{true}$ with $z^\mathrm{true}_i \sim \mathcal{N}(0, 1)$ with probability $0.1$ and otherwise $z^\mathrm{true}_i = 0$.
We choose $\mathcal{X}$ to be the smallest box containing the measurements $D z^\mathrm{true}$.
Both methods start from the minimum-norm least-squares solution $s^0$ that satisfies $(D^\tpose D)s^0 = D^T \bar{x}$, where $\bar x$ is a randomly selected measurement.

\paragraph{GPU code generation.}
The instance of Table~\ref{tab:codegen_ista_fista} differs from the verified one in its size and data.
It has $m = 10^4$, $n = 4 \cdot 10^4$, and $\lambda = 0.02$.
The nonzeros of each row of $D$ lie in positions given by $65$ random permutations of the columns, which gives $260$ entries per row before duplicates are summed and $257$ to $260$ after.
Their values are drawn i.i.d.\ from $\mathcal{N}(0, 1/260)$.
The measurement $x$ ranges over the box $\{x \mid \|x - \bar x\|_\infty \le 1\}$ around a center $\bar x$ drawn from $\mathcal{N}(0, I)$.
Both methods use the step size $\eta = 1/(\|D\|_1 \|D\|_\infty) = 0.0147$ and the initial iterate $s^0 = 2 \ones$, with $200$ iterations per solve in double precision.
Each of seven interleaved timing rounds takes the median of three SciPy solves and the mean of ten solves of the CuPy implementation and of each generated solver, and we report the median over rounds.
Each timed solve of a generated CUDA solver includes uploading the initial iterate and downloading the solution, while the matrices stay on the device.
All rounds ran on one node with two AMD EPYC 7413 CPUs and an NVIDIA A100 SXM4 GPU with $80$\,GB of memory, in a job with $16$ cores and the node to itself.
The OpenMP solvers use $16$ threads, and the SciPy implementation runs on one core.
The generated code is compiled with GCC 11.5 and CUDA 13.3, for the A100 only.
The agent ran on a login node of the cluster, and its builds, checks, and timings ran on GPU nodes with NVIDIA A100 GPUs with $80$\,GB of memory, either SXM4 GPUs with AMD EPYC 7413 CPUs or PCIe GPUs with Intel Xeon CPUs of the Ice Lake generation.
Each solver is generated three times, in independent agent runs, with the default budgets of Appendix~\ref{app:codegen_files}: $10$ minutes for translation and $2$ hours for each of at most five rounds of code optimization.

\subsection{Network optimization}\label{app:experiments_network}

We consider a minimum cost network flow problem with varying demands~\cite{bertsekas1991linear}.
The network has $n_s$ supply nodes and $n_d$ demand nodes, joined by $n_e$ capacitated directed edges from supply to demand nodes.
Each edge $e = (i, j)$ has a gain $g_e > 0$: one unit of flow shipped on it leaves supply node $i$ and delivers $g_e$ units at demand node $j$.
The supply and demand constraint matrices $A_s \in \reals^{n_s \times n_e}$ and $A_d \in \reals^{n_d \times n_e}$ therefore have the entries $(A_s)_{ie} = 1$ and $(A_d)_{je} = -g_e$ for each edge $e = (i, j)$, and zeros elsewhere.
Gains model losses and conversions along an edge, such as yields or exchange rates.
We write $A \in \reals^{(n_s + n_d) \times n_e}$ for the matrix that stacks $A_s$ on top of $A_d$.
When the gains spread over orders of magnitude, $A$ is badly scaled, whereas a plain incidence matrix, with gains equal to one, is already equilibrated.
For each edge $i$, the scalar $c_i$ is the cost of sending one unit of flow and the scalar $\bar z_i$ is the flow capacity.
The vector $b_s \ge 0$ collects the supply available at each supply node, and the vector $x \le 0$ collects the demand requested at each demand node.
With the flow $z \in \reals^{n_e}$ as the optimization variable, the problem is the linear program
\begin{equation}\label{prob:network_flow}
	\begin{array}{ll}
		\text{minimize} & c^\tpose z\\
		\text{subject to} & A_s z \le b_s,\\
		& A_d z = x,\\
		& 0 \le z \le \bar z,
	\end{array}
\end{equation}
which routes the flow at minimum cost while respecting the supply, the demand, and the capacities.
The demand $x$ is the parameter, so each instance of the family asks the same network to serve a different load.

\paragraph{Algorithms.}
We solve~\eqref{prob:network_flow} with the PDHG algorithm~\cite{chambolleFirstOrderPrimalDualAlgorithm2011}, in the form adapted to linear programs~\cite{applegatePracticalLargeScaleLinear2021}.
The algorithm operates on the primal-dual iterate $s^k = (z^k, v^k, y^k)$, where $v^k$ is the dual variable of the supply inequality and $y^k$ is the dual variable of the demand equality.
With a fixed step size $\eta > 0$, the iterations are
\begin{equation}\label{eq:pdhg_iterations}
	\begin{array}{rl}
		z^{k+1} &= \sat_{[0, \bar z]}\!\left(z^k - \eta (c + A_s^\tpose v^k - A_d^\tpose y^k)\right),\\
		v^{k+1} &= \relu\!\left(v^k + \eta \left(A_s (2 z^{k+1} - z^k) - b_s\right)\right),\\
		y^{k+1} &= y^k + \eta \left(x - A_d (2 z^{k+1} - z^k)\right),
	\end{array}
\end{equation}
with step size $\eta = 0.8 / \|A\|_2$, computed from the matrix the method iterates on.
One iteration of~\eqref{eq:pdhg_iterations} composes a saturated linear step, a ReLU step, and an affine step, so it mixes two different PWA primitives.
We also analyze a variant with momentum in the primal variable, which replaces $2 z^{k+1} - z^k$ by $2 \tilde{z}^{k+1} - z^k$ in both dual updates, where
\begin{equation}\label{eq:pdhg_momentum_iterations}
    \tilde{z}^{k+1} = z^{k+1} + \xi_{k+1} (z^{k+1} - z^k),
\end{equation}
with $\xi_{k+1} = k / (k + 3)$, and $\xi_{k+1} = 0$ recovers PDHG.
The extrapolation~\eqref{eq:pdhg_momentum_iterations} is an affine step that \fomtools{} eliminates, so the momentum variant costs the user one extra line and the VP no extra variables.

\paragraph{Equilibration.}
Two-sided Ruiz equilibration~\cite{ruizScalingAlgorithmEquilibrate2001} of $A$, with $30$ passes, returns positive diagonal matrices $E = \diag(E_s, E_d)$ and $D$ such that the rows and columns of $E A D$ have nearly unit $\ell_\infty$-norm.
The equilibrated method runs the iterations~\eqref{eq:pdhg_iterations} on the scaled problem, with the data
\begin{equation*}
    \hat A_s = E_s A_s D, \quad \hat A_d = E_d A_d D, \quad \hat c = D c, \quad \hat{\bar z} = D^{-1} \bar z, \quad \hat b_s = E_s b_s,
\end{equation*}
and with the parameter entering as $E_d x$.
Its iterate $\hat s^k = (\hat z^k, \hat v^k, \hat y^k)$ corresponds to the original iterate $(z^k, v^k, y^k) = (D \hat z^k, E_s \hat v^k, E_d \hat y^k)$.
The parameter set is not scaled, so all variants share the same $\mathcal{X}$ and the same samples.
The matrices $E$ and $D$ depend only on the fixed matrices $A_s$ and $A_d$, so they are constants, and each iteration is still one saturated linear step, one ReLU step, and one affine step.
On the verified instance equilibration reduces $\|A\|_2$ from $9.72$ to $2.54$, which lengthens the step size by the same factor, and reduces the ratio of the largest to the smallest nonzero entry magnitude from $81.4$ to $17.2$.

\paragraph{Residual.}
For PDHG and its momentum variant, the residual map
\begin{equation*}
    \hat H_K = - \begin{bmatrix}
        \eta^{-1} I & -\hat A_s^\tpose & \hat A_d^\tpose\\
        -(1 + 2\xi_K) \hat A_s & \eta^{-1} I & 0\\
        (1 + 2\xi_K) \hat A_d & 0 & \eta^{-1} I
    \end{bmatrix}
\end{equation*}
turns the last step into the residual of the problem the method iterates on.
Writing each update of~\eqref{eq:pdhg_iterations} as an optimality condition, as in~\cite[App.~B.2]{ranjanExactVerificationFirstOrder2025}, gives
\begin{equation*}
    \hat H_K (\hat s^K - \hat s^{K-1}) \in
    \begin{bmatrix}
        \hat c + \hat A_s^\tpose \hat v^K - \hat A_d^\tpose \hat y^K + \mathcal{N}_{[0, \hat{\bar z}]}(\hat z^K)\\
        \hat b_s - \hat A_s \hat z^K + \mathcal{N}_{\reals_+}(\hat v^K)\\
        \hat A_d \hat z^K - E_d x
    \end{bmatrix}.
\end{equation*}
For the equilibrated method this is the residual of the \emph{scaled} problem, whose rows are measured in scaled units, so its norm cannot be compared with that of the unscaled method.
We therefore unscale the residual before it enters the objective, and measure it with
\begin{equation}\label{eq:network_residual}
    \begin{gathered}
    H_K = \diag(D^{-1}, E_s^{-1}, E_d^{-1}) \, \hat H_K,
    \\
    H_K (\hat s^K - \hat s^{K-1}) \in
    \begin{bmatrix}
        c + A_s^\tpose v^K - A_d^\tpose y^K + \mathcal{N}_{[0, \bar z]}(z^K)\\
        b_s - A_s z^K + \mathcal{N}_{\reals_+}(v^K)\\
        A_d z^K - x
    \end{bmatrix}.
    \end{gathered}
\end{equation}
The inclusion follows by substituting the scaled data, since $\hat A_s^\tpose \hat v^K = D A_s^\tpose v^K$ and $\hat A_s \hat z^K = E_s A_s z^K$, and similarly for $A_d$.
Positive diagonal scaling leaves the normal cones of $[0, \bar z]$ and $\reals_+$ unchanged.
The right-hand side of~\eqref{eq:network_residual} is the residual of the original problem~\eqref{prob:network_flow} at the unscaled iterate: the dual residual, the supply slack, and the demand mismatch, all in the original units.
For the unscaled methods $E$ and $D$ are identities and $H_K = \hat H_K$.
All four VPs therefore maximize the same quantity, the original residual, over the same parameter set, so their optimal values are directly comparable.
The unscaling is a constant diagonal matrix applied to the residual expression passed to \texttt{ft.norm}. 

\paragraph{Verified instance.}
The network is a random bipartite graph with $n_s = 10$ supply nodes and $n_d = 15$ demand nodes, in which each of the $150$ possible edges is present independently with probability $0.3$.
The realized graph has $n_e = 49$ edges, and every node has between $2$ and $8$ of them.
The gains are drawn log-uniformly from $[0.1, 10]$, the costs $c_e$ uniformly from $[1, 10]$, and the capacities $\bar z_e$ uniformly from $[2, 6]$, with random seed $285$, and every supply node has supply $6$.
The demand ranges over $\mathcal{X} = [-6, -4]^{n_d}$.
We take the initial iterate to be a primal-dual solution of the instance at the center of the parameter set, $x = -5 \ones$, mapped to the scaled variables for the equilibrated method.
The instance $x = -6 \ones$ has a solution, so every instance in the family is feasible.
We solve the VP for $K = 1, \ldots, 39$.

\paragraph{Loop versus unrolled solvers.}
The solvers of Table~\ref{tab:codegen_network_schedule} differ from the verified instance in the network and the number of iterations.
The network is a complete bipartite graph with $n_s = 20$ supply nodes and $n_d = 30$ demand nodes, so $n_e = 600$, with gains, costs, and capacities drawn as above with random seed $20$, and every supply node has supply $4$.
The demand ranges over $[-6, -4]^{n_d}$, and the initial iterate is a primal-dual solution at its center.
We run momentum PDHG~\eqref{eq:pdhg_momentum_iterations} with Ruiz equilibration and the step size $\eta = 0.5/\|\hat A\|_2 = 0.0927$ for $K = 50$ iterations.
The loop solver's iteration block has $2$ steps and $15$ coefficient arrays with $8{,}520$ nonzeros, together with the table of the $50$ values of $\xi_k$.
The unrolled representation of the full run has $149$ steps and $217$ coefficient arrays with $150{,}230$ nonzeros.
The memory columns count $8$ bytes per table entry and per stored nonzero, and the library column gives the size of the optimized shared library.
Both solvers agree with the NumPy reference execution to within $10^{-14}$, and the timing follows the protocol at the start of this appendix.

\subsection{Optimal control}\label{app:experiments_opt_control}

We consider the control of a constrained time-invariant linear dynamical system over a horizon of $h$ time steps~\cite[Ch.\ 11.3]{borrelliPredictiveControlLinear2017}.
At time $t$, the system has state $v_t \in \reals^{n_v}$ and input $u_t \in \reals^{n_u}$.
Given dynamics matrices $A^\mathrm{dyn}$ and $B^\mathrm{dyn}$, a positive semidefinite state cost $Q$ with terminal cost $Q_h$, and a positive definite input cost $R$, the control problem is
\begin{equation}\label{prob:optimal_control}
	\begin{array}{ll}
		\text{minimize} & (1/2) \left( v_h^\tpose Q_h v_h + \sum_{t=1}^{h-1} \left( v_t^\tpose Q v_t + u_t^\tpose R u_t \right) \right)\\
		\text{subject to} & v_{t+1} = A^\mathrm{dyn} v_t + B^\mathrm{dyn} u_t, \quad t = 1, \ldots, h-1,\\
		& v_t \in \mathcal{V}, \quad t = 1, \ldots, h,\\
		& u_t \in \mathcal{U}, \quad t = 1, \ldots, h-1,\\
		& v_1 = x,
	\end{array}
\end{equation}
where the optimization variables are the state and input sequences, and the sets $\mathcal{V}$ and $\mathcal{U}$ are boxes.
The parameter $x$ is the initial state, which the controller measures afresh at every sampling instant.
Eliminating the states through the dynamics turns~\eqref{prob:optimal_control} into an equivalent problem in the inputs alone, with variable $z = (u_1, \ldots, u_{h-1}) \in \reals^{(h-1) n_u}$~\cite[Ch.\ 11.3]{borrelliPredictiveControlLinear2017},
\begin{equation}\label{prob:control_condensed_qp}
	\begin{array}{ll}
		\text{minimize} & (1/2) z^\tpose P z + q(x)^\tpose z\\
		\text{subject to} & l(x) \le M z \le u(x),
	\end{array}
\end{equation}
where the matrix $P$ is positive definite and the vectors $q(x)$, $l(x)$, and $u(x)$ depend affinely on the initial state.
A model predictive controller solves~\eqref{prob:control_condensed_qp} once per sampling period, so the sampling rate fixes the iteration budget $K$ and the question of \Sec\ref{sec:introduction} becomes a hard real-time requirement.

\paragraph{Algorithms.}
The method of Lu and Yang~\cite{luPracticalOptimalFirstOrder2026} solves programs with one-sided inequality constraints, so we stack the two-sided rows of~\eqref{prob:control_condensed_qp} as $G z \le b(x)$, with $G = [M^\tpose \; {-M^\tpose}]^\tpose$ and $b(x) = (u(x), -l(x))$, and drop the rows whose bound is infinite.
We solve the result with their restarted accelerated primal-dual hybrid gradient method (rAPDHG), which runs inner loops of $K_\mathrm{r}$ iterations.
Each inner loop starts from a restart point $(z^0, y^0)$, with $z^{-1} = z^0$ and running averages $\bar z^0 = z^0$ and $\bar y^0 = y^0$, and at inner index $k = 0, \ldots, K_\mathrm{r} - 1$ it iterates
\begin{equation}\label{eq:rapdhg_iterations}
	\begin{array}{rl}
		z_\mathrm{md}^k &= (1 - \beta_k^{-1}) \bar z^k + \beta_k^{-1} z^k,\\
		y^{k+1} &= \Pi_{\reals_+}\!\left( y^k + \tau_k \left( G \left( z^k + \theta_k (z^k - z^{k-1}) \right) - b(x) \right) \right),\\
		z^{k+1} &= z^k - \eta_k \left( P z_\mathrm{md}^k + q(x) + G^\tpose y^{k+1} \right),\\
		\bar z^{k+1} &= (1 - \beta_k^{-1}) \bar z^k + \beta_k^{-1} z^{k+1}, \quad \bar y^{k+1} = (1 - \beta_k^{-1}) \bar y^k + \beta_k^{-1} y^{k+1},
	\end{array}
\end{equation}
with $\beta_k = (k + 2)/2$, $\theta_k = k/(k + 1)$, $\eta_k = (k + 1) / \left( 2 w (\|P\|_2 + K_\mathrm{r} \|G\|_2) \right)$, and $\tau_k = w (k + 1) / (2 K_\mathrm{r} \|G\|_2)$ for a primal weight $w > 0$.
After $K_\mathrm{r}$ iterations the method restarts from $(\bar z^{K_\mathrm{r}}, \bar y^{K_\mathrm{r}})$.
One iteration uses a single ReLU step, on the rows of $G$, and every other update is affine, the restart included, so the iteration solves no linear system.

The projection onto a closed convex set $C$ satisfies $w - \Pi_C(w) \in \mathcal{N}_C(\Pi_C(w))$~\cite[Proposition 16.44]{bauschkeConvexAnalysisMonotone2017}.
Applying this property to the dual update and rearranging the primal update gives the scaled residual, as in the derivations of~\cite{ranjanExactVerificationFirstOrder2025},
\begin{equation}\label{eq:rapdhg_residual}
	\begin{bmatrix}
		\eta_k^{-1} (z^k - z^{k+1}) - P (z_\mathrm{md}^k - z^{k+1}) \\
		\tau_k^{-1} (y^k - y^{k+1}) - G \left( (z^{k+1} - z^k) - \theta_k (z^k - z^{k-1}) \right)
	\end{bmatrix}
	\in
	\begin{bmatrix}
		P z^{k+1} + q(x) + G^\tpose y^{k+1} \\
		b(x) - G z^{k+1} + \mathcal{N}_{\reals_+}(y^{k+1})
	\end{bmatrix},
\end{equation}
where $\mathcal{N}_{\reals_+}$ is the normal cone of the nonnegative orthant and the first block holds with equality.
A zero residual certifies a KKT point of the condensed program, and dividing by the step sizes removes the factor $k + 1$ by which the steps grow within an inner loop.
We measure it at the inner iterate of every step, restart steps included, and the user obtains it by stacking the two expressions in one call to \texttt{ft.norm}, as in \Sec\ref{subsec:objectives}.

\paragraph{Verified instance.}
We use the linearized dynamics of the quadcopter of~\cite{kouzoupisProperAssessmentQP2015}, as given in the OSQP documentation~\cite{osqpDocs2026}, sampled every $0.1$ seconds, with $n_v = 12$ states, $n_u = 4$ inputs, and a horizon of $h = 5$ time steps, so the condensed program has $16$ variables.
The state and input cost matrices are $Q = \diag(0, 0, 10, 10, 10, 10, 0, 0, 0, 5, 5, 5)$, $Q_h = 10 Q$, and $R = 0.1 I$.
The inputs are the rotor thrusts relative to hover, with $\mathcal{U} = \{ u \in \reals^4 \mid -0.9916 \ones \le u \le 2.4084 \ones \}$, and the set $\mathcal{V} = \{ v \in \reals^{12} \mid |v_1| \le \pi/6, \; |v_2| \le \pi/6 \}$ bounds the two tilt angles and leaves the other states free.
The initial state lies in
\begin{equation*}
    \mathcal{X} = \{ x \in \reals^{12} \mid |x_1| \le \pi/9, \; |x_2| \le \pi/9, \; |x_6| \le 1/2, \; x_i = 0, \; i \ne 1, 2, 6 \},
\end{equation*}
\ie, hovering positions with tilt angles of at most $20$ degrees and an altitude within $0.5$ meters of the hover reference.
We compare rAPDHG with $K_\mathrm{r} = 6$ and $K_\mathrm{r} = 12$, both with primal weight $w = 1$ and started from $z^0 = 0$ and $y^0 = 0$.
We solve the VP for $K = 1, \ldots, 75$.

\paragraph{Closed-loop flight.}
The controller of Figure~\ref{fig:control_codegen_flight} differs from the verified instance in its sampling period, horizon, constraints, and restart period.
It uses the same quadcopter, cost matrices $Q$, $Q_h = 10 Q$, and $R$, and input box $\mathcal{U}$, with the dynamics sampled every $0.01$ seconds and a horizon of $100$ steps.
We obtain these dynamics from the continuous-time model recovered by a matrix logarithm of the $0.1$-second dynamics, followed by zero-order-hold discretization.
We eliminate the states as in~\eqref{prob:control_condensed_qp} and hold the input constant over blocks of $10$ consecutive steps, so the variable $z$ has $40$ entries.
The cost tracks a reference trajectory, so the vector $q$ depends affinely on both the measured state and the reference over the horizon.
We bound the two tilt angles by $\pi/6$ at the last state of each input block, which gives $40$ rows of $G$, and with the $80$ rows of the input box the matrix $G$ has $120$ rows and the vector $b(x)$ depends on the measured state.
We run rAPDHG~\eqref{eq:rapdhg_iterations} with $K_\mathrm{r} = 64$ and $w = 1$.

The reference is the figure eight $(\sin \omega t, \; (1/2) \sin 2 \omega t, \; 1 + 0.4 \sin(\omega t / 2))$, in meters, with $\omega = 2\pi/4$, so one horizontal lap takes $4$ seconds and the full reference repeats every $8$ seconds.
Each flight lasts $8$ seconds and starts on the reference, and each implementation runs the number of iterations that fits in the budget at its measured time per iteration.
The simulated plant applies the first input clipped to $\mathcal{U}$, and the tilt constraints hold only at the ends of the input blocks.

\subsection{Portfolio rebalancing}\label{app:experiments_portfolio}

An investor holds the portfolio $x \in \reals^n$, whose entries are the fractions of one unit of capital held in each of $n$ assets, and chooses new holdings $h \in \reals^n$.
Given a positive definite covariance matrix $\Sigma$, a forecast $\mu$ of the asset returns, a risk aversion $\kappa > 0$, a transaction cost $\lambda > 0$, and a vector $\bar z > 0$ of position limits, the rebalancing problem is
\begin{equation}\label{prob:portfolio}
	\begin{array}{ll}
		\text{minimize} & (\kappa/2) h^\tpose \Sigma h - \mu^\tpose h + \lambda \|h - x\|_1\\
		\text{subject to} & \ones^\tpose h = 1,\\
		& 0 \le h \le \bar z,
	\end{array}
\end{equation}
where $h$ is the optimization variable, and the current holdings $x$ and the return forecast $\mu$ are the parameters.
The risk aversion $\kappa$ trades the variance of the portfolio against its expected return, and the last term charges a proportional cost $\lambda$ on each trade $h_i - x_i$.
The covariance follows a factor model $\Sigma = FF^\tpose + dI$, with a factor loading matrix $F \in \reals^{n \times k}$ and $d > 0$.
As in~\cite{schallerEmbeddedCodeGeneration2022}, we write the risk as $h^\tpose \Sigma h = \|f\|_2^2 + d\|h\|_2^2$, with the $k$ factor exposures $f = F^\tpose h$ as auxiliary variables.

The parameter set contains the holdings that one period of bounded market moves can produce from a reference portfolio, which we take to be the balanced portfolio $x^\mathrm{b} = (1/n)\ones$.
Suppose the investor holds $x^\mathrm{b}$ at the start of a period and does not trade during it.
Asset $i$ earns the return $r_i$, the relative change in its price, so the capital held in asset $i$ grows from $1/n$ to $(1 + r_i)/n$.
The total capital becomes $\sum_j (1 + r_j)/n$, so the fraction of it held in asset $i$ is
\[
x_i = \frac{1 + r_i}{\sum_{j=1}^n (1 + r_j)}, \qquad i = 1, \ldots, n.
\]
We call $x$ a drifted portfolio, since it differs from $x^\mathrm{b}$ only because prices moved.
It is still fully invested, $\ones^\tpose x = 1$, but it holds more than $1/n$ in the assets whose return exceeded the average and less in the others, and the rebalancing problem~\eqref{prob:portfolio} decides how to trade it back.
For example, with $n = 2$ assets and returns $r = (0.1, -0.1)$, the drifted portfolio is $x = (0.55, 0.45)$.
We assume that every return lies in $[-\bar r, \bar r]$ for a given $\bar r \in [0, 1)$.
Coordinate $x_i$ increases with $r_i$ and decreases with every other return, so it is largest when asset $i$ gains $\bar r$ and every other asset loses $\bar r$, and smallest in the opposite case.
Every drifted portfolio therefore lies in the polyhedron
\begin{gather*}
    \mathcal{X} = \{ x \in \reals^n \mid \underline{x} \ones \le x \le \overline{x} \ones, \; \ones^\tpose x = 1 \},\\
    \underline{x} = \frac{1 - \bar r}{(1 - \bar r) + (n - 1)(1 + \bar r)}, \qquad
    \overline{x} = \frac{1 + \bar r}{(1 + \bar r) + (n - 1)(1 - \bar r)},
\end{gather*}
the intersection of the smallest box containing them with the budget hyperplane.
The forecast has the form $\mu_i = \bar\mu_i + 0.1\, \Sigma_{ii}^{1/2} s_i$, where $\bar\mu$ is a fixed mean forecast and the score $s_i$ lies in $[-3, 3]$.
The parameters $(x, \mu)$ therefore range over $\mathcal{X} \times \mathcal{B}$, where $\mathcal{B} = \{ \mu \in \reals^n \mid |\mu_i - \bar\mu_i| \le 0.3\, \Sigma_{ii}^{1/2}, \; i = 1, \ldots, n \}$.

We write the transaction cost with an auxiliary variable $\tau \in \reals^n$ that satisfies $\tau \ge |h - x|$ elementwise, and replace $\lambda \|h - x\|_1$ by $\lambda \ones^\tpose \tau$.
With the stacked variable $z = (h, \tau, f) \in \reals^{2n + k}$, problem~\eqref{prob:portfolio} takes the form
\begin{equation}\label{prob:condensed_qp}
	\begin{array}{ll}
		\text{minimize} & (1/2) z^\tpose P z + q(\mu)^\tpose z\\
		\text{subject to} & l(x) \le M z \le u(x),
	\end{array}
\end{equation}
where the matrix $P$ is positive semidefinite and the vectors $q(\mu)$, $l(x)$, and $u(x)$ depend affinely on the parameters.
For the portfolio problem these data are
\begin{gather*}
    P = \kappa \begin{bmatrix} d I & 0 & 0\\ 0 & 0 & 0\\ 0 & 0 & I \end{bmatrix}, \quad
    q(\mu) = \begin{bmatrix} -\mu\\ \lambda \ones\\ 0 \end{bmatrix},
    \\
    M = \begin{bmatrix} \ones^\tpose & 0 & 0\\ F^\tpose & 0 & -I\\ I & 0 & 0\\ I & -I & 0\\ I & I & 0 \end{bmatrix}, \quad
    l(x) = \begin{bmatrix} 1\\ 0\\ 0\\ -\infty\\ x \end{bmatrix}, \quad
    u(x) = \begin{bmatrix} 1\\ 0\\ \bar z\\ x\\ +\infty \end{bmatrix}.
\end{gather*}
The five groups of rows of $M$ are the budget, the factor exposures $F^\tpose h - f = 0$, the position limits, the sell rows $h - \tau \le x$, and the buy rows $h + \tau \ge x$, and the last two together encode $\tau \ge |h - x|$.
The holdings enter only the one-sided bounds of the sell and buy rows, and the forecast enters only the linear cost, so $P$ and $M$ are fixed.

\paragraph{Algorithms.}
We solve~\eqref{prob:condensed_qp} with the ADMM scheme of the OSQP solver~\cite{stellatoOSQPOperatorSplitting2020}.
Given parameters $\rho, \sigma > 0$, the algorithm iterates on $s^k = (z^k, v^k)$ through
\begin{equation}\label{eq:osqp_iterations}
	\begin{array}{rl}
		w^{k+1} &= \sat_{[l(x), u(x)]}(v^k),\\
		(P + \sigma I + \rho M^\tpose M) z^{k+1} &= \sigma z^k - q(\mu) + \rho M^\tpose (2 w^{k+1} - v^k),\\
		v^{k+1} &= v^k + M z^{k+1} - w^{k+1}.
	\end{array}
\end{equation}
Setting $y^k = \rho(v^{k-1} - w^k)$ makes the pair $(w^k, y^k)$ satisfy the normal cone conditions of~\eqref{prob:condensed_qp} at every iteration~\cite{banjacInfeasibilityDetectionAlternating2019,stellatoOSQPOperatorSplitting2020}, so the only optimality conditions left to measure are primal and dual feasibility. 
These can be written as a linear map of the fixed-point residual~\cite{ranjanExactVerificationFirstOrder2025},
\begin{equation}\label{eq:osqp_residual_scaling}
    \begin{bmatrix} P z^K + M^\tpose y^K + q(\mu)\\ M z^K - w^K \end{bmatrix}
    = H (s^K - s^{K-1}),
    \qquad
    H = \begin{bmatrix} -\sigma I & -\rho M^\tpose\\ 0 & I \end{bmatrix}.
\end{equation}
Over-relaxation with $\alpha \in (0, 2)$ replaces $z^{k+1}$ in the second line by the solution $\tilde z$ of the linear system, and then updates
\begin{equation*}
	\begin{array}{rl}
		z^{k+1} &= \alpha \tilde z + (1 - \alpha) z^k,\\
		v^{k+1} &= v^k + \alpha (M \tilde z - w^{k+1}).
	\end{array}
\end{equation*}
The value $\alpha = 1$ recovers~\eqref{eq:osqp_iterations}, for which the identity~\eqref{eq:osqp_residual_scaling} holds.

\paragraph{Problem data.}
The instances have $n = 200$, $500$, $1000$, or $3000$ assets and $k = 100$ factors, with $\kappa = 1$, $\lambda = 0.01$, $\bar z = 0.3 \ones$, $d = 0.1$, and the return bound $\bar r = 0.15$.
The entries of $F$ are i.i.d.\ normal with mean zero and variance $1/k$, and the entries of $\bar\mu$ are drawn uniformly from $[0, 0.5]$, both with NumPy's default generator seeded with $(0, n, k)$.
Each solve draws a drifted portfolio $x$ with returns uniform on $[-\bar r, \bar r]$, and a forecast $\mu$ whose scores $s_i$ are standard normal and truncated to $[-3, 3]$.
In the CVXPY model, the risk uses the factor exposures, and $x$ and $\mu$ are parameters whose boxes $\underline{x} \ones \le x \le \overline{x} \ones$ and $\mathcal{B}$ are passed through \texttt{param\_bounds}.
CVXPYgen converts the problem to a quadratic program with $2n + k$ variables and $4n + k + 1$ constraint rows, in which the parameters enter only the linear cost and the upper bounds of the inequality rows.

Both solvers run the over-relaxed form of the iterations~\eqref{eq:osqp_iterations} on the Ruiz-equilibrated problem, with $\rho = 1$, $\sigma = 10^{-6}$, over-relaxation $1.6$, and $\rho$ scaled by $10^3$ on the equality rows, starting from zero at every solve.
CVXPYgen generates OSQP with $\rho = 0.1$ and an adaptive choice of $\rho$, and it exposes neither setting.
We therefore set $\rho = 1$ with \texttt{osqp\_update\_rho} before the first solve and turn the adaptation off.
Every solver stops once the $\ell_\infty$-norms of the primal and dual residuals of the unscaled problem are at most $10^{-4}$, checked every $25$ iterations, with at most $4000$ iterations.
For OSQP, this is its test with absolute tolerance $10^{-4}$ and relative tolerance zero.
By default, OSQP~1.x also requires the duality gap to fall below the absolute tolerance, and we turn this check off.
The \fomtools{} solvers stop on the primal and dual residuals~\eqref{eq:residuals} of the KKT conditions, which the verifier also uses as its performance metric.
OSQP instead measures the primal residual at its auxiliary variable, as $\|Mz - w\|_\infty$, which is never smaller since $w \in [l(x), u(x)]$.
All solves reach the tolerance.
The mean iteration counts are $99$, $109$, $144$, and $291$ for OSQP and $100$, $110$, $158$, and $298$ for both \fomtools{} solvers, for $n = 200$, $500$, $1000$, and $3000$.
The agent ran with budgets of $10$ minutes for translation and $20$ minutes for each round of code optimization, the second below the default of $2$ hours in Appendix~\ref{app:codegen_files}.
Each solver is generated three times, in independent agent runs.
Generating the translated solver from the CVXPY problem takes $41$ to $323$ seconds, and adding the code optimization phase, in one to five rounds, brings the total to $8$ to $43$ minutes.

Each timing round runs every solver through a C driver, which sets all entries of $x$ and $\mu$ through CVXPYgen's update functions and then calls its solve function.
The driver makes one warm-up pass over a sequence of $100$ draws and then times three passes.
We report the median over three interleaved rounds of the mean time per solve, on an Apple M4 Pro, and then the median over the three agent runs.
The C times stay within $2.5\%$ of their medians across rounds.

\end{document}